\documentclass[10pt]{amsart}
\usepackage{amssymb, amscd, amsmath, amsthm, epsf, epsfig, latexsym}

\def\zz{{\bf Z}}
\def\ff{{\bf F}} \def\hh{{\bf H}}
\def\qq{{\bf Q}}
\def\cc{{\bf C}}
\def\rr{{\bf R}}
\def\co{\colon\thinspace}

\def\calm{\mathcal{M}}
\def\calg{\mathcal{G}}
\def\co{\colon}

\newcommand\Aut{\operatorname{Aut}}
\newcommand\Hom{\operatorname{Hom}}
\newcommand{\fig}[2] { \includegraphics[scale=#1]{#2} }

\newtheorem{theorem}{Theorem}[section]
\newtheorem{lemma}[theorem]{Lemma}
\newtheorem{corollary}[theorem]{Corollary}
\newtheorem{prop}[theorem]{Proposition}
\theoremstyle{definition}
\newtheorem{definition}[theorem]{Definition}

\numberwithin{equation}{section}

\begin{document}

\title[Geography of 4--manifolds]{The Geography Problem for 4--Manifolds with Specified Fundamental Group}

\author{Paul Kirk}
\author{Charles Livingston}
\thanks{This work was supported by grants from the NSF}
\address{Department of Mathematics, Indiana University, Bloomington, IN 47405}
\email{pkirk@indiana.edu}
\email{livingst@indiana.edu}




\begin{abstract} 
For any class $\calm$ of 4--manifolds, for instance the class $\calm(G)$ of closed oriented manifolds with $\pi_1(M) \cong G$ for a fixed group $G$, the geography of $\calm$ is the set of integer pairs $\{(\sigma(M), \chi(M))\ | \ M \in \calm\}$, where $\sigma$ and $\chi$ denote the signature and Euler characteristic. This paper explores general properties of the geography of $\calm(G)$ and undertakes an extended study of $\calm(\zz^n)$.
\end{abstract}

\maketitle

\section{Introduction}

The fundamental group of a 4--manifold constrains its other algebraic invariants in interesting ways.  The results of this article concern  the  constraints imposed by the fundamental group of a 4-manifold $M$ on its  two fundamental numerical invariants, its Euler characteristic $\chi(M)$ and signature $\sigma(M)$.

The starting point of our investigation into this topic was  Hausmann and Weinberger's~\cite{hw} construction of a perfect group  $G$ for which any 4--manifold $M$ with $\pi_1(M) \cong G$ satisfied $\chi(M) > 2$.  (As a consequence  it follows that Kervaire's~\cite{ke} classification of groups of homology $n$--spheres $\Sigma^n,   n   \ge 5$ does not extend to dimension four.) Work extending that of  Hausmann and Weinberger, investigating large classes of groups, applying the techniques of $L^2$-homology, and studying the problem of finding constraints   on $\chi(M)$ and  $\sigma(M)$  arising from $\pi_1(M)$, include~\cite{e1,e2,jk,ko,  luck}.

Aspects of this study have appeared in different guises, especially when focused on particular categories of manifolds.    For instance, for symplectic manifolds we have the initial result of Gompf~\cite{go}  that any finitely presented group arises as the fundamental group of a symplectic manifold, followed by work  of Baldridge and Kirk~\cite{Baldridge-Kirk} refining Gompf's construction to build manifolds with relatively small second Betti number.

From a different perspective, the study of 4-manifolds can be  divided into two parts. The first  part consists of identifying all possible homotopy types of 4--dimensional Poincar\'e complexes with given fundamental group.  The second part consists of identifying for each such Poincar\'e complex $X$,  the set of 4--manifolds homotopy equivalent to $X$.   As an example, for the trivial group,  Whitehead's theorem 
\cite{jhc} states that for each unimodular   symmetric form there is a simply connected 4--dimensional Poincare complex $X$ with that intersection form and $X$ is unique up to homotopy equivalence.  Freedman's work~\cite{fq}  implies that   for each $X$ there are either one or two closed topological manifolds homotopy equivalent to $X$, depending on whether the form is even or odd, respectively.  Rochlin's theorem implies that in the smooth category there are some $X$ which cannot be realized by a smooth 4--manifold. Donaldson's work~\cite{d, dk} placed further constraints on the forms that are realizable  and demonstrated nonuniqueness; extending that work, it has been shown  that for some forms there are an infinite number of non-diffeomorphic manifolds realizing the form~\cite{fm, ov}.    When the group is nontrivial, the problem becomes much more difficult, even in the topological setting; see~\cite{hkt} for the results in the case $G \cong\zz/2$.  Beyond extremely simple groups, the possibility of classifying manifolds with a given fundamental group is completely intractable.  Studying the signatures and Euler characteristics of such manifolds is a more accessible first step in which a number of interesting examples and important problems arise.  

\medskip
\vskip.1in
\noindent{\bf Summary of main results.}

The  study of the Euler characteristics and signatures of 4--manifolds with specified group can be encompassed by the general study of the geography of a group.  We begin with some basic definitions;   assume throughout that all spaces and maps are based, manifolds are closed and oriented,  and that groups $G$ are finitely presented.  For now manifolds can be taken in the topological category, though our results hold in other settings as well.

\begin{definition} \label{def4} Fix a  group $G$ and $\alpha \in H_4(G)$.
 \begin{enumerate} 
 \item $\calm(G)$ denotes the class of  4--manifolds $M$ with fundamental group isomorphic to $G$. 
 \item $\calm(G,\alpha)$ denotes the class of pairs $(M, f)$ where $f \co \pi_1(M) \to G$ is an isomorphism and  $f_*([M])=\alpha $. 
 \end{enumerate}
  \end{definition}

We associate to $\calm(G)$ and $\calm(G,\alpha)$ subsets of $\zz^2$:

\begin{definition} \label{def5} Fix a finitely presented group $G$ and class $\alpha\in H_4(G)$.
 \begin{enumerate} 
 \item $\calg(G) = \{(\sigma(M), \chi(M)) \ |\ M \in \calm(G) \}.$ 
 \item $\calg(G, \alpha) = \{(\sigma(M), \chi(M)) \ |\ (M,f)  \in \calm(G, \alpha) \}.$ 
 \end{enumerate} 
 \end{definition}

We will note later that for distinct $f, g \co \pi_1(M) \to G$,   $f_{*}([M])= \Psi_* g_{*}([M])$ where $\Psi_*$ is the automorphism of $H_4(G)$ induced by the automorphism   $\Psi = f \circ g^{-1}$ of $G$, and that for any such automorphism $\Psi$ of $G$, $\calg(G, \alpha) = \calg(G, \Psi_*(\alpha))$.

The {\em geography problem for $G$} is the problem of identifying the sets $\calg(G)$. The similarly defined geography problem for $\calg(G, \alpha)$ is a  natural extension; observe  that  $\calg(G, \alpha)\subset \calg(G)$ and  $\calg(G)=\cup_{\alpha\in H_4(G)}\calg(G,\alpha).$    The structure of the geography of a group will be seen to be captured by the following  functions.

\begin{definition} $ \ $
\begin{enumerate}
\item $q_G(\sigma)$:   the minimum value of $\chi$ for pairs $(\sigma, \chi) \in \calg(G)$.

\item $q_{G,\alpha}(\sigma)$:  the minimum value of $\chi$ for pairs $(\sigma, \chi) \in \calg(G, \alpha)$.
\end{enumerate}
\end{definition}

 In the first part of this paper  we prove a number of basic results concerning the geography problem.   Some of this work constitutes a survey of past work, usually restated in terms of the general geography problem, while some of the work is new, especially the development of properties of $\calg(G, \alpha)$.  In the second part, we investigate the problem in the special classes of groups $G$, in particular the case of $G$ free abelian, building on work in our earlier article~\cite{kl}.  

Sections \ref{secnot} and \ref{secgeo} present general results concerning the geography problem. We begin by explaining why $\calm(G)$ and $\calm(G,\alpha)$ are non-empty and prove the following (for a precise statement see Theorem \ref{qthm}), a natural generalization of past work on the subject. 

\medskip

\noindent{\bf Theorem. } {\em For any finitely presented group $G$, the function $q_{G}(\sigma)$ completely determines $\calg(G)$, the Hausmann-Weinberger invariant $q(G)$, the Kotschick invariant $p(G)$ and the function $f_G(t)$ of \cite{Baldridge-Kirk}.  Moreover, for any $\alpha\in H_4(G)$,}
$$q_{G,\alpha}(\sigma)\ge q_G( \sigma) \ge \max\{|\sigma| -2\beta_1(G) +2, \beta_2(G)-2\beta_1(G)+2\}.$$

\noindent In addition, various symmetry properties of $q_{G}$ and $q_{G,\alpha}$ are established in  Section \ref{secgeo}.  As mentioned earlier, for an automorphism $\Psi$ of $G$, there is an induced automorphim $\Psi_*$ of $H_4(G)$.  We explain in Section \ref{secgeo} the identification of $\calg(G, \alpha) $ and $\calg(G, \Psi_*(\alpha))$.

One of our goals is  establishing sharper bounds than provided by Theorem \ref{qthm}.  Our approach  is based on 
identifying the part of the intersection form of $M\in \calm(G)$ determined entirely by the group $G$. More precisely, if $f\co \pi_1(M)\to  G$ is an isomorphism, let $I(M,f)\subset H^2(M)$ denote the image of $f^* \co H^2(G)\to H^2(M)$. Also, 
for a symmetric bilinear form $\phi$, let $b^{\pm}(\phi)$ denote the maximal dimension of a  subspace on which $\phi$ is positive (respectively, negative) definite; let $b^{\pm}(M)$ 
denote the corresponding numbers for the intersection form of $M$.

Corollary \ref{cor1} implies the following theorem, which generalizes the fact that the intersection form of an aspherical 4-manifold is determined by its fundamental group. 

\medskip

\noindent{\bf Theorem.}  {\em Let $\alpha\in H_4(G)$ and let $(M,f) \in \calm(G, \alpha)$. Then the restriction of the intersection form of $M$ to $I(M,f)$ depends only on $G$ and $\alpha\in H_4(G)$; it is equivalent to the pairing   $\phi(x,y) = (x\cup y)\cap \alpha$ on $H^2(G)$. In particular, 
 $b^\pm(M)\ge b^\pm(\phi)$ and any isotropic subspace of the pairing   $\phi$   is also isotropic for the intersection form of $M$. If $\phi$ is even, then the restriction of the intersection form of $M$ to $I(M,f)$ is also even.  }

\medskip

In Section \ref{secex} we turn to the problem of calculating $q_G(\sigma)$ and $q_{G,\alpha}(\sigma)$ for various classes of groups $G$, in part illustrating the extent to which Corollary~\ref{cor1} can be used to improve previously known bounds.  Furthermore, the explicit calculation of  $q_G(\sigma)$ (and $q_{G,\alpha}(\sigma)$) calls on the construction of  manifolds realizing the given bounds.   We prove some general results in the case when $\alpha=0$ and when $\alpha$ is a multiple.   Theorems \ref{surface} and \ref{3manifolds}  imply the following.
\medskip

\noindent{\bf Theorem.}  {\em 
\begin{enumerate}
\item  If $G$ is the fundamental group of a closed orientable surface $F$ of genus greater than 0, then $q_G(\sigma)=|\sigma| + 2\chi(F)$.  

\item If $G$ is the fundamental group of a closed oriented 3--manifold, and $G= F_n*H$, the free product of a free group  $F_n$ on $n$ generators and $H$ is not a free product of any group with a free group, then $q_{G}(\sigma)=|\sigma|+ 2-2n$. 
\end{enumerate} }

\medskip

Note that for 2-- and 3--manifold groups, $H_4(G)=0$, so that in this theorem $\alpha$ does not arise.

We then turn to a detailed investigation of free abelian groups, $G \cong \zz^n$.  Early work on the subject did not extend successfully beyond $n = 4$, with the exception of our work in~\cite{kl} which determined the Hausmann-Weinberger function   $q(\zz^n)$ for all $n$.  Recall the    definition:
$$q(G) = \min \{ \chi(M) \ |\  M \in  \calm (G)\} = {\min}_\sigma q_G(\sigma).$$
(We will show later that   $q(\zz^n) = q_{\zz^n}(0)$.)

 The cohomology algebra of $\zz^n$ is an exterior algebra, and hence the bilinear  form $\phi$ described above can be interpreted as follows. Fix an identification $\Lambda^n(\zz^n)\cong \zz$. Then given $\omega\in \Lambda^{n-4}(\zz^n)$ (the Poincar\'e dual of $\alpha$) one obtains a bilinear pairing:
 $$\Lambda^2(\zz^n)\times \Lambda^2(\zz^n)\to \zz, \ (x,y)\mapsto x\wedge y\wedge \omega.$$
Despite this simple formula, the properties of this pairing (such as its rank, determinant, and signature) are difficult to extract. We achieve some success for $n\leq 6$ by studying the $\Aut(\zz^n)$ action on $\Lambda^*(\zz^n)$ to put the class $\omega$ in standard form.   This gives us extensive information about  $\calg(\zz^n,\alpha)$ for $n\leq 6$ and $\alpha\in H_4(\zz^n)$. The calculations are recorded in Theorems
\ref{nis4}, \ref{nis5}, and \ref{nequals6}. The  most detailed example of this article  concerns the geography of the group $\zz^6$.  A complete analysis of $q_{\zz^6}(\sigma)$ is achieved, illustrating phenomena that have not been observed before.    We  also carry out a detailed analysis for some classes $\alpha \in H_4(\zz^6)$, revealing further the complexity of the geography problem, even for fairly simple groups, and underscoring a number of outstanding questions.

In Section \ref{asym}, we investigate  the asymptotics of $\calg(\zz^n)$ as $n\to \infty$.   One consequence of  our earlier work \cite{kl} is that $q(\zz^n)/n^2\to \frac{1}{2}$ as $n\to \infty$.  Kotschick's invariant $p(\zz^n)$ equals $\min_\sigma (q_{\zz^n}(\sigma)-\sigma)$.  Constructing manifolds $M$ with $\pi_1(M) \cong \zz^n$ for which the signature is large relative to the second Betti number has proved to be extremely difficult, and there are no general results concerning the problem.  An initial conjecture might be that such manifolds do not exist, and that  $p(\zz^n)/n^2\to \frac{1}{2}$.  However we are able to produce counterexamples:     Theorem \ref{13/28thm}  states the following.

\medskip

\noindent{\bf Theorem.} {\em For an infinite set of $n$, $p(\zz^n)/n^2 \le  \frac{13}{28}$.}
 
\medskip

The paper concludes with Section~\ref{problemsec} organizing some of the outstanding problems related to the geography of groups.    

\medskip

\noindent{\it Acknowledgments.}  The authors would like to thank P.~Teichner, I.~Hambleton, and C.~Van Cott  for many helpful comments, and M.~Larsen for providing   insights into our study of the asymptotics related to  Kotschick's function $p(G)$ in Section~\ref{13/28thm}.  


\section{Notation and Basic Results}\label{secnot}

We have thus far been vague about what category of 4--manifold one considers. There are many possible choices, such as smooth 4--manifolds, topological 4--manifolds, symplectic 4--manifolds, compact K\"ahler surfaces, irreducible or minimal 4--manifolds;  each choice leads to an interesting set of questions. Moreover, delicate questions arise when comparing the problems for different categories.  One can also leave the realm of manifolds and consider  4--dimensional Poincar\'e complexes.  Most of our discussion applies in this setting, and we anticipate that to fully understand the case of topological 4--manifolds, the case of 4--dimensional Poincar\'e complexes will need to be analyzed  so that surgery theory~\cite{fq} might be applied.

The focus of the present article is on topological and smooth manifolds. Usually, the manifolds we construct to realize a pair $(\sigma,\chi)\in \calg(G,\alpha)$ will be smooth, and our results about which pairs $(\sigma,\chi)$ cannot be realized will apply to 4--dimensional Poincar\'e complexes.

The geography problem is often stated in terms of other characteristic numbers: in complex geometry one considers the {\em  holomorphic Euler characteristic} $\chi_h(M)$  and $c_1^2(M)$. Since $\chi_h(M)=\tfrac{1}{4}(\sigma(M)+\chi(M))$ and $c_1^2(M)= 3\sigma(M) + 2\chi(M)$ whenever these are defined, the pairs $(\chi_h, c_1^2)$ and $(\sigma, \chi)$ contain the same information. For our purposes the latter formulation is more convenient.

\subsection{Euler characteristic}

For a space $X$ we denote the integral and rational homology by $H_*(X)$ and $H_*(X; \qq)$, respectively. The Betti numbers are defined by $\beta_i (X)= \mbox{dim} (H_i(X; \qq))$. If the homology is finitely generated,   the {\em Euler characteristic} is given by $\chi(X) = \sum (-1)^i \beta_i(X)$. For a closed, orientable 4--manifold $M$, Poincar\'e duality implies that 
\begin{equation}\label{euler1} \chi(M)=2-2\beta_1(M)+ \beta_2(M). \end{equation}

\subsection{Symmetric bilinear forms and the signature}

For convenience we recall some terminology and facts about symmetric bilinear integer forms.

We will be considering symmetric bilinear forms $F \co V \times V \to \zz$ with $V$ a finitely generated free abelian group, whose rank is called the {\em rank of }$F$. If the determinant of $F$ (in some basis) is $\pm 1$ the form is called {\em unimodular}; if the determinant is non-zero the form is called {\em non-degenerate}. If $F(v,v)$ is even for all $v\in V$, the form is called {\em even}, and otherwise it is called {\em odd}. The {\em signature} of $F$ is the difference \begin{equation}\sigma(F)=b^+(F)-b^-(F), \end{equation} where $b^+(F)$ (respectively $b^-(F)$) is the dimension of the largest subspace of $V\otimes \rr$ on which (the obvious extension of) $F$  is positive (respectively negative) definite. A non-degenerate form is called {\em definite} if $|\sigma(F)|$ equals the rank of $F$, and called {\em indefinite} otherwise.

Two basic unimodular (integral) forms are the indefinite even form of rank two and signature 0, denoted $H$ (represented by a $2 \times 2$ matrix with diagonal entries 0 and off-diagonal entries 1) and the even positive definite form of rank 8 (and signature 8), denoted $E_8$.

The classification theorem (see for example~\cite{MH}) states that any unimodular indefinite odd form is equivalent to a diagonal form with all diagonal entries $\pm 1$, and thus is determined by its rank and signature. Any unimodular indefinite even form is a direct sum $kH \oplus mE_8$, where $k$ is a positive integer, and $m$ is an integer, and thus is again determined by its rank and signature. The classification of definite forms is not complete; all such even forms have rank and signature divisible by 8.

An {\em isotropic} subspace of a symmetric bilinear form $F\co V\times V\to \zz$ is a subspace $W$ so that $F(w_1,w_2)=0$ for all $w_1,w_2\in W$. If $F$ is non-degenerate and admits an $n$--dimensional isotropic subspace, an easy argument shows that $b^+(F)\ge n$ and $b^-(F)\ge n$, so that $$\text{rank}(F)\ge |\sigma(F)|+2n.$$

\subsection{Some algebraic topology}

For a finitely presented group $G$ there is an associated Eilenberg-MacLane space which we denote $B_G$. By definition, $B_G$ is a based, connected CW-complex satisfying $\pi_1(B_G) \cong G$ and $\pi_i(B_G) = 0, i \ge 2$; $B_G$ is unique up to homotopy. The homology groups $H_*(G)$ and the cohomology ring $H^*(G)$ are defined to be $H_*(B_G)$ and $H^*(B_G)$, where coefficients are assumed to be integers unless specifically noted.

Every compact $m$--manifold has the homotopy type of a finite CW-complex. (In fact, by results of Kirby-Siebenmann~\cite{KS}, every compact manifold $M$ has the homotopy type of a finite polyhedron.) Thus, given a homomorphism $\phi\co \pi_1(M) \to G$, there is a corresponding based map, unique up to based homotopy, $f\co M \to B_G$ inducing $\phi$ on fundamental groups. This proves the following well-known  lemma.

\begin{lemma} Let $M$ be a closed, oriented $m$--manifold and $G$ a discrete group. A homomorphism $\phi\co\pi_1(M)\to G$ uniquely determines an element $\alpha=f_*([M])\in H_m(G)$, where the  map $f\co M\to B_G$ is the unique based homotopy class inducing $\phi$ and $[M]\in H_m(M)$ denotes the orientation class.\qed

\end{lemma}

Thus we will refer to $\alpha\in H_m(G)$ as the {\em class induced by} $\phi\co \pi_1(M)\to G$.

Bilinear parings can be constructed on the cohomology $H^*(G)$ using the cup and cap products as follows. Given a class $\alpha\in H_{2n}(G)$, define \begin{equation}\label{eq1} H^{n}(G)\times H^n(G)\to \zz \text{ by } (x,y)\mapsto (x\cup y)\cap \alpha.\end{equation} When $n$ is even, the pairing (\ref{eq1}) is symmetric.

\begin{definition} Suppose that $M$ is an oriented closed $2n$--dimensional manifold and $\phi \co \pi_1(M)\to G$ is a homomorphism, with $f \co M\to B_G$ the corresponding homotopy class of maps, and $\alpha=f_*([M])\in H_{2n}(G)$ the class induced by $\phi$. 

Define the subspace 
$$\hskip-.6in  I(M,f) = \text{ image }( f^*\co H^{n}(G)\to H^{n}(M)\to  H^{n}(M)/\text{torsion}) $$
$$\hskip2.5in \subset H^n(M)/\text{torsion}.$$

\end{definition}

The bilinear pairing (\ref{eq1}) completely determines the restriction of the intersection form of $M$, $$H^n(M)/\text{torsion}\times H^n(M)/\text{torsion}\to \zz,\ (x,y)\mapsto (x\cup y)\cap [M],$$ to $I(M,f)$, since \begin{equation}\label{eq5}(f^*(a) \cup f^*(b)) \cap [M] = (a \cup b) \cap \alpha. \end{equation} In particular, if $f^* \co H^{n}(G)\to H^{n}(M)$ is surjective, then the intersection form of $M$ is completely determined by algebra, or, more precisely, by the cohomology ring $H^*(G)$.

The problem of deciding when a class $\alpha\in H_{m}(G)$ is represented by a map $f\co M\to B_G$ from a closed $m$--manifold is a classical (and difficult) problem in topology (see for example~\cite{thom}), but can always be solved for $m=4$. We outline the argument.

\begin{lemma} Given any group $G$ and class $\alpha\in H_4(G)$, there exists an oriented closed smooth 4--manifold $M$ and a continuous map $f\co M\to B_G$ so that $f_*([M])=\alpha$. \end{lemma}

\begin{proof} The proof of this follows from a calculation using  the bordism spectral sequence,  which is itself  an application of the Atiyah-Hirzebruch spectral sequence (see~\cite{st} or~\cite{gw})  applied to the generalized homology theory given by oriented bordism, $\Omega_*^{SO}$.   In greater detail, there is a spectral    sequence, with $E_2$-term given by \{$H_i(X,\Omega_j^{SO})\}$, converging to $\Omega_*^{SO}(X)$.  The relevant coefficient groups, oriented bordism groups of a point, are given by $\Omega_0^{SO} = \zz$,  $\Omega_1^{SO} = 0$, $\Omega_2^{SO} = 0$, $\Omega_3^{SO} = 0$, and $\Omega_0^{SO} = \zz$, generated by ${\bf CP^2}$ (see~\cite{ms}).  Thus, it follows from the spectral sequence that there is an exact sequence $ \Omega_4^{SO} \to \Omega^{SO}_4(G) \to H_4(G,\zz)\to 0$, giving the desired surjection.

\end{proof}

The problem of determining the size of $I(M,f)$ is also difficult in general, but the following fact gives a simple criterion when $M$ is a 4--manifold.

\begin{theorem}\label{thm1} Let $\alpha\in H_4(G)$ and suppose $(M,f)\in \calm(G,\alpha)$. Then the homomorphism $f_* \co H_{2}(M) \to H_2(G)$ is surjective and $f^* \co H^2(G) \to H^2(M)$ is injective. In particular, rank$(I(M,f))=\beta_2(G)$. \end{theorem}

\begin{proof} An Eilenberg-MacLane space for $G$ can be built from $M$ by adding cells of dimension $3$ and higher. It follows that the inclusion $H_2(M) \to H_2(B_G)$ is surjective. The result for cohomology follows from the Universal Coefficient Theorem. \end{proof}

Surgery yields a method to ensure that a given class $\alpha \in H_m(G)$ is represented by a map $f \co M^m\to B_G$ {\it inducing an isomorphism on fundamental groups}.  The following argument  is roughly that of Wall~\cite[Theorem 1.2]{wa}, for the most part translated into combinatorial group theory.

\begin{lemma} 
Suppose $G$ is finitely presented and $\alpha \in H_m(G)$, $m>3$, is represented by 
$f \co M \to B_G$ for a closed oriented $m$--manifold $M$. Then it is also represented by $g \co N^m \to B_G$ so that the induced morphism $g_* \co \pi_1(N) \to G$ is an isomorphism. 
\end{lemma}

\begin{proof} 
We first show that if $\phi:H\to G$ is a homomorphism of finitely presented groups, then $\phi$ can be extended to an epimorphism $\Phi:H*F\to G$  (where $H*F$ is the free product of $H$ with a finitely generated free group $F$) so that the kernel of $\Phi$ is normally generated by finitely many elements. 

Indeed, given presentations 
$$H=\langle h_1,\cdots, h_k\ | \ w_1, \cdots, w_\ell\rangle
\text{ 
 and } G=\langle g_1,\cdots , g_m\ | \ r_1,\cdots , r_n \rangle,$$
let $F$ denote the free group generated by $g_1,\cdots , g_m$. Then 
$$H*F=\langle h_1,\cdots , h_k, g_1,\cdots, g_m\ | \ w_1,\cdots, w_\ell\rangle$$
and there is an obvious epimorphism of $H*F$ to $G$, taking $h_i$ to $\phi(h_i)$ and $g_i \in F$ to $g_i \in G$. 
For $i=1,\cdots, k$, let $z_i\in F$ be a word in the $g_i \in F$ which is sent to $\phi(h_i) \in G$ by the canonical surjection $F\to G$. 
Taking the quotient of $H*F$ by the subgroup generated by the finitely many elements
 $r_1,\cdots, r_n,   z_1h_1^{-1},\cdots,  z_kh_k^{-1}$ clearly yields $G$, since the generators  $h_i$ and relations $ z_jh_j^{-1}$ can be eliminated.

Thus by replacing $M$ by the connected sum of $M$ with finitely many copies of $S^1\times S^{n-1}$ we may arrange that $\pi_1(M)\to G$ is surjective and has kernel normally generated by finitely many elements. This can be done without changing the homology class $\alpha$,  by arranging that the maps of $S^1\times S^{n-1}$ to $B_G$ factor  through $S^1$.

Since $M$ is orientable and has dimension greater than three, the set of normal generators of the kernel can be represented by disjointly embedded circles with trivial normal bundles, and surgery on these circles yields a manifold $N$ with $\pi_1(N)\cong G$. Extending $M\to B_G$ to the trace of the surgery shows that $N$ maps to $B_G$, inducing an isomorphism on fundamental groups, without altering the class $\alpha\in H_m(G)$. \end{proof}

\begin{corollary} $\calm(G,\alpha)$ is non-empty for any $ \alpha\in H_4(G)$.\qed \end{corollary}

\begin{corollary} \label{cor1}Let $\alpha\in H_4(G)$ and let $(M,f)\in \calm(G,\alpha)$. Then the restriction of the intersection form of $M$ to $I(M,f)$ depends only on $G$ and $\alpha\in H_4(G)$; it is given by the pairing (\ref{eq1}):   $(x,y) \mapsto (x\cup y)\cap \alpha$. Moreover, 

\begin{enumerate} \item $b^\pm(M)\ge b^\pm(I(M,f))$, \item any isotropic subspace of the pairing (\ref{eq1})    is also isotropic for the intersection form of $M$, and \item if (\ref{eq1}) is even, then the restriction of the intersection form of $M$ to $I(M,f)$ is also even. \end{enumerate}

Thus if $I(M,f)=H^2(M)/\text{torsion}$, (for instance, if   $f^* \co H^2(G)\to H^2(M)$   is surjective), then the intersection form of $M$ is determined by $\alpha$. \qed

\end{corollary}

Corollary \ref{cor1} forms the starting point of our investigation of the sets $\calg(G)$ and $\calg(G,\alpha)$. It should be thought of as a generalization of the fact that the intersection form of an aspherical 4--manifold is determined by its fundamental group.


\section{Geography and fundamental groups} \label{secgeo}
\subsection{Basic properties of $\calg(G)$} We now apply the observations of the previous section to the problem of the geography associated to the class of 4--manifolds with a specified fundamental group.

First, note that Equation (\ref{euler1}) and the fact that $\beta_1(X)=\beta_1(G)$ for a connected space $X$ satisfying $\pi_1(X)=G$ implies that if $M$ is a closed, orientable 4--manifolds $M$ with fundamental group $G$, then

\begin{equation}\label{euler2} \chi(M) = 2 - 2 \beta_1(G) + \beta_2(M). \end{equation}

We have the following basic observations. \begin{theorem}\label{basicthm} If $M$ is a closed 4--manifold and $\pi_1(M) \cong G$, then \begin{enumerate} \item $ \chi(M)\equiv \sigma(M)\pmod{2}.$ \item $ \chi(M) \ge \beta_2(G) - 2\beta_1(G) + 2 .$ \item $ \chi(M) \ge |\sigma(M)| - 2\beta_1(G) + 2.$ \end{enumerate} \end{theorem}

\begin{proof} We have that $\beta_2(M) = b^+(M) + b^-(M)$ and $\sigma(M) = b^+(M) - b^-(M)$. The mod 2 congruence  then follows from Equation (\ref{euler2}).

Theorem \ref{thm1} shows that $\beta_2(M)\ge \beta_2(G)$. Together with Equation (\ref{euler2}) this proves the first inequality.

Again, using that $\beta_2(M) = b^+(M) + b^-(M)$ and $\sigma(M) = b^+(M) - b^- (M)$, we have that $$\chi(M) = 2 - 2 \beta_1(G) \mp \sigma(M) +2b^\pm(M).$$ The second inequality now follows from the fact that $b^\pm(M) \ge 0$. \end{proof}

Theorem \ref{basicthm} gives us our first interesting constraint on the geography.

\begin{corollary}\label{basiccor} Fix $G$ and $\alpha\in H_4(G)$. Then $$\calg(G,\alpha)\subset \calg(G)$$ and $\calg(G)$ is a subset of the intersection of the three half-planes $$\calg(G)\subset \{ \chi\ge \beta_2(G) - 2\beta_1(G) + 2 \}\cap \{ \chi+\sigma\ge - 2\beta_1(G) + 2 \}\cap\{ \chi-\sigma\ge - 2\beta_1(G) + 2 \}.$$
 
Moreover, if $(\sigma,\chi)\in \calg(G,\alpha)$, then $(\sigma+1,\chi+1), (\sigma-1,\chi+1), (\sigma,\chi+2)\in \calg(G,\alpha).$ \end{corollary}

\begin{proof} All but the last assertion follows from Therorem \ref{basicthm}. For the last assertion, observe that a map $f\co M\to B_G$ can be homotoped to be constant on a 4--ball in $M$. Construct a map on the connected sum $f'\co M\# \cc P^2\to B_G$ by making $f'$ constant on the $\cc P^2$ factor. Since $\cc P^2$ is simply connected, $f$ induces an isomorphism on fundamental groups. Moreover, $f'$ and $f$ determine the same class $\alpha\in H_4(M)$ since $f'$ is constant on $\cc P^2$. Hence if $f\co M\to B_G$ represents $(\sigma,\chi)\in \calg(G,\alpha)$, $f'\co M\#\cc P^2\to B_G$ represents $(\sigma+1,\chi+1)\in \calg(G,\alpha)$. Similar arguments using $-{\cc P}^2$ and $S^2\times S^2$ show that $(\sigma-1,\chi+1), (\sigma,\chi+2)\in \calg(G,\alpha).$

\end{proof}

\begin{figure}[h]
\centerline{\fig{.4}{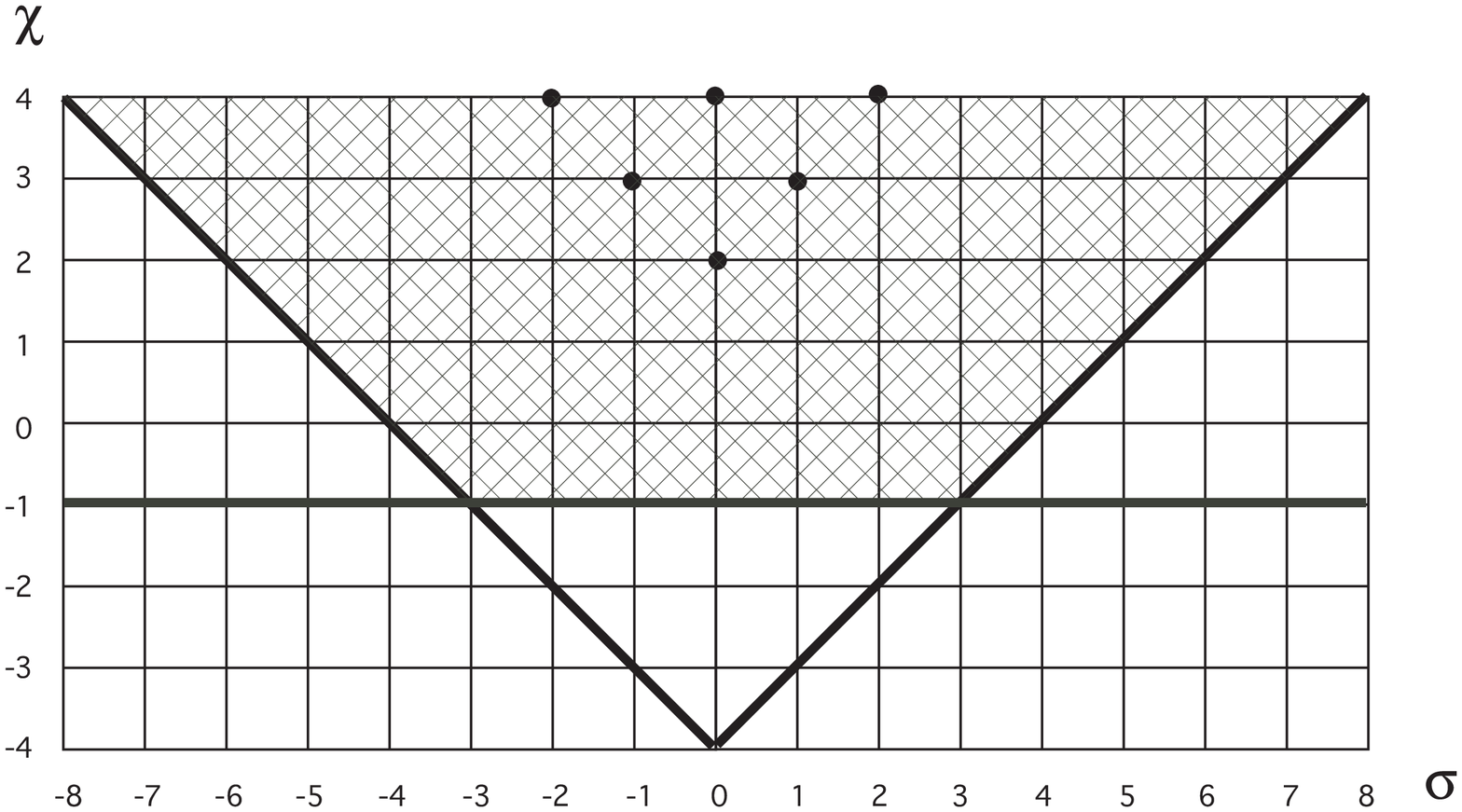}}
\caption{The bounds of Corollary~\ref{basiccor} and the geography of  $G \protect\cong \zz^3$}\label{boundfig}
\end{figure}

We will show below that these bounds are not strict in the case of $G \cong \zz^n, n \ge 2$.  As an example, in Figure~\ref{boundfig} we illustrate with  three dark lines the bounds given by Corollary~\ref{basiccor} for $G \cong \zz^3$ while the black dots indicate the only values that   actually occur in the region, as will be shown later.

\begin{prop} \label{prop4.3}Let $\Psi \co G\to G$ be an automorphism of $G$, inducing an automorphism $\Psi_*\co H_*(G)\to H_*(G)$. Then $\calg(G,\alpha)=\calg(G,\Psi_*(\alpha))$. In other words, $\calg(G,\alpha)$ depends only on the orbit of $\alpha$ under the action of  $\Aut(G)$ on $H_4(G)$. 
Moreover, the fundamental class of any oriented closed 4--manifold $M$ with $\pi_1(M)$ isomorphic to $G$  determines a well-defined element in the orbit set $H_4(G)/\Aut(G)$.   
\end{prop}
\begin{proof} Any automorphism $\Psi\co G\to G$ is induced by a based homotopy equivalence $h:B_G\to B_G$. Hence if $(M,f)\in \calm(G,\alpha)$, then $(M,h\circ f)\in \calm(G,\Psi_*(\alpha))$ and so $\calg(G,\alpha)\subset \calg(G,\Psi_*(\alpha))$.  The reverse inclusion follows by considering $\Psi^{-1}$. 

Given a manifold $M$ with $\pi_1(M)$ isomorphic to $G$, a choice of isomorphism determines a map $f\co M\to B_G$ inducing the given isomorphism on fundamental groups. Another choice gives a possibly different map $g\co M\to B_G$.  The automorphism $ \Psi=g_*\circ f_*^{-1}:\pi_1(B_G)\to \pi_1(B_G)$ then shows that $\Psi_*(f_*([M]))=g_*([M])$. Hence $f_*([M])$ is well-defined in $H_4(G)/\Aut(G)$.  
\end{proof}

A consequence of Proposition \ref{prop4.3} is that if we denote by $[\alpha]\in H_4(G)/\Aut(G)$  the orbit of $\alpha$, then it makes sense to define $\calm(G,[\alpha])$ to be the class of 4--manifolds  $M$ whose fundamental group is isomorphic to $G$ such that $f_*([M])\in [\alpha]$ for some (and hence any) reference map $f:M\to B_G$. In particular, $\calm(G)$ is the {\em disjoint} union of the $\calm(G,[\alpha])$, and the geography $\calg(M,[\alpha])$ is well-defined.  Put another way, there is a well--defined surjective assignment $\calm(G)\to H_4(G)/\Aut(G)$; this takes homotopy equivalent $M_1,M_2\in \calm(G)$ to the same $[\alpha]\in H_4(G)/\Aut(G)$.  

A simple bordism argument pointed out to the authors by Peter Teichner shows that   two smooth manifolds $M_1$ and $M_2$ are in $\calm(G,[\alpha])$ if and only if     the connected sum of $M_1$ with sufficiently many copies of $\cc P^2$ and $-\cc P^2 $ is diffeomorphic to  the connected sum of $M_2$ with sufficiently many copies of $\cc P^2$ and $-\cc P^2 $.

\medskip

Reversing orientation changes the sign of $\sigma$ and preserves $\chi$, and hence one obtains the following.

\begin{prop}\label{symm} For any $G$, $\calg(G)$ is symmetric  with respect to reflection
 through the $\chi$--axis. More generally, $\calg(G,-\alpha)$ is obtained by reflecting $\calg(G,\alpha)$ 
through the $\chi$--axis.  Moreover, if $\alpha\in H_4(G)$ and $G$ admits an automorphism $\Psi$ so that $\Psi_*(\alpha)=-\alpha$, then $\calg(G,\alpha)$ is symmetric with respect to the $\chi$--axis. \qed\end{prop}

There is no reason why $\calg(G,\alpha)$ should be symmetric for general $\alpha$. (The symmetry assertions in Proposition \ref{symm} extend to geography problems for classes of 4--manifolds closed under change of orientation; this excludes symplectic manifolds, for example.)

\subsection{The invariants of Hausmann-Weinberger and Kotschick}

We now explain the relationship between the problem of identifying $\calg(G)$ with previously studied invariants, notably the Hausmann-Weinberger invariant~\cite{hw} $$q(G)=\inf \{ \chi(M) \ | \ M\in \calm(G) \}$$ and its variants due to Kotschick~\cite{ko} $$p(G)=\inf \{ \chi(M)- \sigma(M) \ | \ M\in \calm(G)\}$$ and Baldridge-Kirk~\cite{Baldridge-Kirk} $$f_G(t)= \inf \{ \chi(M)+ t \sigma(M) \ | \ M\in \calm(G)\}, t\in \rr.$$ Note first that $q(G)=f_G(0)$ and $p(G)=f_G(-1)$.

The Hausmann-Weinberger invariant can be viewed in the present context as the smallest $\chi$--intercept of any horizontal line  which intersects $\calg(G)$ nontrivially. Kotschick's variant is similar, but considers lines of slope $1$ rather than slope 0.  More generally, $f_G(t)$ is the smallest $\chi$--intercept of a line of slope $-t$ that intersects $\calg(G)$ nontrivially.

Corollary \ref{basiccor} shows that translating by $(0,2)$ preserves $\calg(G)$ and $\calg(G,\alpha)$. This motivates the introduction of the following functions.

\begin{definition} $\ $
\begin{enumerate}
\item  $q_G(\sigma) = \min\{\chi(M)\ | \ M \in \calm(G) \text{\ and\ } \sigma(M) = \sigma\}.$ \item $q_{G,\alpha}(\sigma) = \min\{\chi(M)\ | \ (M,f) \in \calm(G,\alpha) \text{\ and\ } \sigma(M) = \sigma\}.$
\end{enumerate}
 \end{definition}

The basic properties of $q(G),p(G),f_G(t),q_G(\sigma)$ and $q_{G,\alpha}(\sigma)$ are summarized in the following theorem.

\begin{theorem}\label{qthm}For each finitely presented group $G$ and $\alpha\in H_4(G)$: \begin{enumerate}
 \item $q_G( \sigma)$ is defined for all $\sigma\in\zz$. \item  $q_G( \sigma) \equiv \sigma \mod 2.$ \item $q_G( \sigma) \ge \max\{|\sigma| -2\beta_1(G) +2, \beta_2(G)-2\beta_1(G)+2\}.$ \item $q_G( \sigma +1) = q_G(\sigma) \pm 1.$

\item For all large values of $\sigma$, $q_G( \sigma) = \sigma +p(G)$  where $p(G)\in \zz$ denotes Kotschick's invariant.

\item The Hausmann-Weinberger invariant $q(G)$ equals $\min\{q_G(\sigma)\ | \ \sigma\in\zz \}$. \item $f_G(t)\ne-\infty$ if and only if  $t\in [-1,1]$. The convex hull of $\calg(G)$ in $\rr^2$ is the intersection of the half planes  $\{\chi+t\sigma\ge f_G(t)\ | \ t\in[-1,1]\}.$  Thus $f_G(t)$ determines and is determined by the convex hull of $\calg(G)$. \item $\calg(G)=\{ (\chi,\sigma)\in \zz^2\  | \  \chi\ge q_G(\sigma) \text{ and } \chi\equiv\sigma\pmod{2}\}$. Thus $q_G(\sigma)$ determines and is determined by $\calg(G)$.

\item $q_G(\sigma)\leq q_{G,\alpha}(\sigma)$.

\item $q_G(-\sigma)=q_G(\sigma)$   and $q_{G,-\alpha}(\sigma)=q_{G,\alpha}(- \sigma).$ \end{enumerate}

The first 4 assertions hold for $q_{G,\alpha}$ also. Moreover,

\begin{enumerate} \item[(5)$^\prime$] There exists integers $p_{+,\alpha}(G)$ and $p_{-,\alpha}(G)$(depending on $G$ and $\alpha$) satisfying $|p_{\pm,\alpha}(G)|\ge p(G)$ such that for all large values of $\sigma$, $q_{G,\alpha}( \sigma) = \sigma +p_{+,\alpha}(G)$ and $q_{G,\alpha}(-\sigma) = -\sigma +p_{-,\alpha}(G)$. \end{enumerate} \end{theorem} \begin{proof} \hfill \begin{enumerate}

\item Given a manifold $M$, by forming the connected  sum with copies of $\pm \cc P^2$ we can build a manifold with the same fundamental group but with $\sigma$ arbitrary. By mapping the extra $\pm\cc P^2$ via the constant map to $B_G$, we can arrange that the class $\alpha\in H_4(G)$ is unchanged. Thus $q_G(\sigma)$ and $q_{G,\alpha}(\sigma)$ are defined for all $\sigma$. \item The parity statement follows from Theorem~\ref{basicthm}. \item This is a restatement of Theorem~\ref{basicthm}. \item Since we can add $\pm \cc P ^2$ to $M$ we see that $|q_G(\sigma +1) - q_G(\sigma)| = 1$.

\item The third assertion shows that the function $\sigma\mapsto q_G(\sigma)- \sigma$ is bounded below as $\sigma\to\infty$. The greatest lower bound  is the Kotschick invariant $p(G)$. The fourth assertion shows that this integer-valued function is non-increasing; in fact, it changes by $0$ or $-2$ when $\sigma$ is replaced by $\sigma+1$. Hence it is eventually constant. The assertion for $q_G$ now follows. For $q_{G,\alpha}$ a separate (but duplicate) argument is needed for $\sigma\to -\infty$ since $q_{G,\alpha}$ may not be symmetric about the $\chi$ axis.

\item This follows from the definition of $q(G)$.

\item Note that $f_G(t)\ne -\infty$ if and only if $t\in [-1,1]$ since Corollary \ref{basiccor} shows that if $(\sigma,\chi)\in \calg(G)$, then so are $(\sigma- n, \chi+n)$ and $(\sigma+n, \chi+n)$ for all $n\in \zz$. For $t\in[-1,1]$, the half plane $\chi+ t\sigma\ge c$ contains $\calg(G)$ whenever $c\leq f_G(t)$. Hence the convex hull of $\calg(G)$ is the intersection of the half planes $\chi + t\sigma\ge f_G(t), \ t\in [-1,1]$.

\item This follows from the definition of $q_G$ as a minimum, the parity statement, and the fact that if $(\sigma,\chi)\in \calg(G)$ then also $(\sigma,\chi+2)\in \calg(G)$. \item Since $\calg(G,\alpha)\subset \calg(G)$, $q_{G,\alpha}(\sigma)\ge q_G(\sigma)$. \item This is seen by reversing orientations.

\end{enumerate} \end{proof}


\section{Basic Examples}\label{secex}

\subsection{The trivial group and free groups} We consider the trivial group, $\{e\}$ and observe that $\#_a\cc P^2\#_b  ({-\cc P}^2) $ realizes any pair $(\sigma, \chi)$ satisfying Corollary~\ref{basiccor} when $a=(\chi+\sigma- 2)/2$ and $b=(\chi-\sigma-2)/2$. In the figure below, we sketch the points of $\calg(\{e\})$ contained in the region $[0,8] \times [0,8]$. We will explain later why the group $\zz^3$ has the same geography. Note that for both of these groups $H_4(G)=0$, so that $\alpha=0$ is the only possibility. Thus $$q_{\{e\}}(\sigma)=|\sigma | + 2=q_{\zz^3}(\sigma).$$

\begin{figure}[h]
\centerline{\fig{.4}{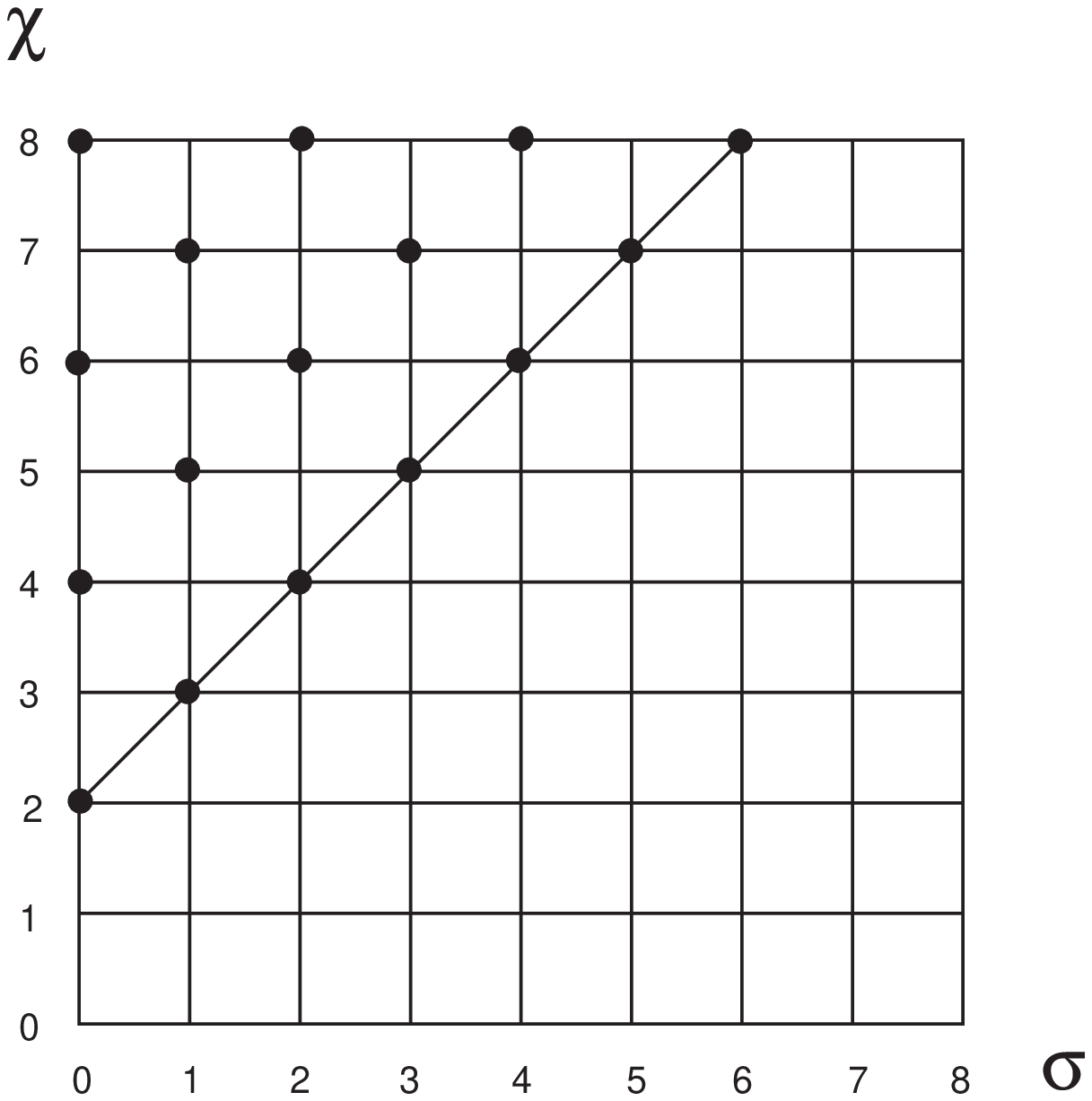}}
\caption{The geography of $\{e\}$  and $\zz^3$}\label{q=0}
\end{figure}

Consider next the free group on $n$ generators, $F_n$. Note that  $H_4(F_n)=0$, $\beta_1(F_n)=n$, and $\beta_2(F_n)=0$. Theorem \ref{qthm} implies that $q_{F_n}(\sigma)\ge |\sigma| - 2n + 2$. The connected sum of $n$ copies of $S^1\times S^3$ with copies of $\pm \cc P^2$ shows that this inequality is an equality: $$q_{F_n}(\sigma)=|\sigma|+2 -2n.$$

\subsection{Torsion classes; 2-- and 3--manifold groups} The {\em deficiency} of a group $G$, def$(G)$, is the infimum over all presentations of $G$ of $ g-r $ , where $g$ is the number of generators and $r$ is the number of relations.

In the case that $H_4(G) = 0$ it is well-known that in the second statement of Theorem~\ref{basicthm},  $\beta_2(G)$ can be replaced with $2\beta_2(G)$:  $\chi(M) \ge 2\beta_2(G) - 2\beta_1(G) +2$.  A proof of a stronger result generalizing to $H_4(G;F)$ with $F$ an  arbitrary field appears in~\cite[Lemma~1.5]{klt}. See Section~\ref{multiplesection} below.  Note that to control the signature it is necessary to work with $\zz$ (or $\qq$) coefficients. .  The idea of the proof is a straightforward generalization  of the earlier proofs. 

\begin{theorem} \label{thm7} Suppose $\alpha\in H_4(G)$ is a torsion class.  Then for any $(M,f)\in \calm(G,\alpha)$, \begin{equation*} q_{G,\alpha}(\sigma)\ge |\sigma| + 2 - 2\beta_1(G) + 2\beta_2(G). \end{equation*}
In particular, if a group $G$ satisfies $H_4(G;\qq)=0$, then this inequality holds for any $\alpha\in H_4(G)$.

If moreover $\alpha=0$, then \begin{equation*} |\sigma| + 2 - 2(\beta_1(G) -\beta_2(G))\leq q_{G,0}(\sigma)\leq |\sigma| +2-2 \ \text{\rm def}(G). \end{equation*} \end{theorem}

\begin{proof}

The cap product $\cap \alpha \co H^4(G)\to \zz$ is zero when $\alpha$ is a torsion class. From Corollary \ref{cor1} and Equation (\ref{eq5}) one concludes that $I(M,f)\subset H^2(M)/\text{torsion}$ is isotropic.

Since the intersection form on $H^2(M)/\text{torsion}$ is unimodular, it follows that $$\text{ rank }(H^2(M))\ge 2\text{ rank }(I(M,f))+|\sigma(M)|.$$ Theorem \ref{thm1} implies that $I(M,f)$ is isomorphic to $H^2(G)/\text{torsion}$, and so the rank of $I(M,f)$ equals $\beta_2(G)$. Thus $\beta_2(M)\ge 2\beta_2(G) +|\sigma(M)|$, which implies the first inequality.

A presentation of a group $G$ with $g$ generators and $r$ relations determines a 2--handlebody $W$ with one 0--handle, $g$ 1--handles, and $r$ 2--handles. Take a map $W\to B_G$ inducing an isomorphism on fundamental groups. Restricting the composite $W\times [0,1]\to W\to B_G$ to $M=\partial (W\times [0,1])$ yields a pair $(M,f)\in \calm(G,0)$ with $\chi(M)=2-2\ \text{def}(G)$ and $\sigma(M)=0$.  Taking connected sums with $\pm \cc P^2$ gives the  upper bound $q_{G,0}(\sigma)\leq |\sigma|+2-2\text{ def}(G)$. \end{proof}

As a quick application of Theorem \ref{thm7}, suppose $G$ is the fundamental group of a closed orientable surface $F$ of genus $g$.  Then $H_4(G)=0$ and hence Theorem \ref{thm7} implies that when $g>0$ (so that $\beta_2(F)=\beta_2(G)$), $$|\sigma| + 4-4g=  |\sigma| + 2 - 2\beta_1(G) + 2\beta_2(G)\leq q_{G}(\sigma)\leq |\sigma| +2-2\ \text{\rm def}(G)=|\sigma|+ 4-4g.$$ This proves the following.

\begin{theorem} \label{surface} If $G$ is the fundamental group of a closed orientable surface $F$ of genus greater than zero, then $H_4(G)=0$ and $$q_G(\sigma)= |\sigma|+ 2\chi(F).$$ \qed \end{theorem}

A more interesting class of examples is given by 3--manifold groups. If $G=\pi_1(N)$ for a compact 3--manifold $N$, then $H_4(G)=0$. Thus Theorem \ref{thm7} applies.

The following theorem completely characterizes the geography of closed oriented 3--manifold groups.

\begin{theorem}\label{3manifolds} Suppose $G$ is the fundamental group of a closed oriented 3--manifold, and assume that $G$ is  a free product $G= F_n*H$, where $F_n$ is a free group on $n$ generators and $H$ is not a free product of any group with non-trivial free group.  Then $$q_G(\sigma)=|\sigma|+2-2n.$$ \qed \end{theorem}

 \begin{proof}  By Kneser's Conjecture~\cite{Hempel}, there is a a closed oriented 3--manifold  $N$   with fundamental group $H$. Performing surgery on the circle $x\times S^1\subset N\times S^1$ yields a 4--manifold $M$ with $\pi_1(M)=H$, $\chi(M)=2$, and $\sigma(M)=0$. Taking the connected sum of $M$ with $n$ copies of $S^1\times S^3$ yields a 4--manifold  with fundamental group $G$, Euler characteristic $2-2n$, and signature zero. Thus $q_G(0)\leq 2-2n$. Taking connected sums with $\pm \cc P^2$ shows that $$q_{G}(\sigma)\leq |\sigma| +2-2n .$$

To prove the reverse inequality, we observe that  the Eilenberg-MacLane space of a free product of two groups is the wedge of the Eilenberg-MacLane spaces of the factors.  Using the Mayer-Vietoris sequence this implies \begin{equation}\label{eq15} \beta_1(G)=\beta_1(F_n) + \beta_1(H)=n + \beta_1(H)\text{ and } \beta_2(G)=\beta_2(F_n) + \beta_2(H) =\beta_2(H). \end{equation}

The prime decomposition and sphere theorems for 3--manifolds imply that $$N=A_1 \#\cdots \# A_k\#K_1 \#\cdots \# K_\ell \# (S^1\times S^2)_1\#\cdots \#(S^1\times S^2)_m$$ where the $A_i$ are aspherical  and the $K_i$ have finite fundamental group  (see for example~\cite{Hempel}). Since we assumed that $H$ is not a free product, the prime decomposition of $N$ cannot have any $S^1\times S^2$ factors and so $$N=A_1 \#\cdots \# A_k\#K_1 \#\cdots \# K_\ell.$$

Since $A_i$ is a closed, orientable, aspherical 3--manifold $$\beta_1(\pi_1(A_i))=\beta_1( A_i) =\beta_2(A_i)=\beta_2(\pi_1(A_i)). $$ Since $K_i$ has finite fundamental group $$ 0=\beta_1(\pi_1(K_i))=\beta_1( K_i) =\beta_2(K_i)\ge\beta_2(\pi_1(K_i))\ge 0. $$ It follows  that $\beta_1(\pi_1(K_i))=0=\beta_2(\pi_1(K_i)).$ Using the Mayer-Vietoris sequence again, we compute \begin{equation}\label{eq16} \beta_1(H)=\sum_j \beta_1(\pi_1(A_i)) =\sum_j \beta_2(\pi_1(A_i)) =\beta_2(H). \end{equation} Equations (\ref{eq15}) and (\ref{eq16})  and Theorem \ref{thm7} imply that $$q_G(\sigma)\ge |\sigma| + 2 - 2n.$$
\end{proof}

 Theorem \ref{3manifolds} does not characterize  closed oriented 3-manifold groups.  One can construct examples as follows. Let  $G$  be a  finite group  with periodic cohomology of period 4.   Then $H_4(G)=0$.  Corollary~4.4 of  the article \cite{hk-4}  by Hambleton and Kreck states that there exists a   topological 4--dimensional rational homology sphere $M$ with fundamental group $G$. Since $\chi(M)=2$ and $\sigma(M)=0$, this example and Theorem \ref{qthm} imply that  
 $q_{G}(\sigma)=2+|\sigma| $.  Milnor   \cite{Msn} listed groups with period 4 cohomology and showed that some of these groups, for example $G$ the symmetric group on three letters,  are not the fundamental groups of 3--manifolds. 
 
Some of these groups, though not fundamental groups of 3--manifolds,  are the fundamental groups of 3--dimensional Poincar\'e complexes.
For instance, Swan has shown in~\cite{swan}  that the symmetric group on three letters is the fundamental group of a 3--dimensional Poincar\'e complex.  Thus these examples do not contradict the possibility that 
  Theorem \ref{3manifolds} characterizes 3--manifold groups in the Poincar\'e category.      In addition, the manifolds constructed in~\cite{hk-4} are topological 4--manifolds and it is not known whether smooth examples exist (see Problem 4.121 of \cite{ki}). Thus an interesting problem is the following.
  
  \medskip
  
  \noindent{\bf Problem.} Find a finitely presented  group $G$  which is not the fundamental group of any Poincar\'e 3--complex (finite or infinite), $G$ is not the free product with a free group,  $H_4(G) = 0$, and   $q_{G}(\sigma)=2+|\sigma| $.    Note that the quesion might have a different answer depennding on whether one considers smooth or topological 4--manifolds. 

\medskip

One can carry out   calculations similar to those of  Theorem \ref{3manifolds} for the fundamental groups of not necessarily closed compact 
3--manifolds; the statements become a bit more complicated. However, one interesting and simple example is the following. 

\begin{theorem} Let $K\subset S^3$ be a knot, and let $G=\pi_1(S^3-K)$. Then $H_4(G)=0$ and
$$q_G(\sigma)=|\sigma|.$$
\end{theorem}
\begin{proof} A knot $K$ in $S^3$ has $H_1(S^3-K)\cong \zz$ and $H_2(S^3-K)=0$. Hence $\beta_1(G)-\beta_2(G)=1$.  The Wirtinger presentation of $G$ has deficiency   1; this  easily implies  def$(G)=1$. The Sphere Theorem \cite{Hempel}   implies that   $S^3-K$ is an Eilenberg-MacLane space, and hence $H_4(G)=H_4(S^3-K)=0$.  The last assertion of Theorem \ref{thm7} then shows that $q_G(\sigma)=|\sigma|$.
\end{proof}

\subsection{Multiple classes}\label{multiplesection}

The following is a well-known result, presented in print in~\cite[Lemma~1.5]{klt}; we state the theorem in terms of classes being multiples, $\alpha = p \tau$, rather than $\alpha = 0 \in H_4(G; \zz / p)$ to be consistent with our approach.  The simple proof is included for the reader's convenience.  

 \begin{theorem}\label{multiple}  Suppose $(M,f) \in\calm(G,\alpha)$ and that $\alpha\in H_4(G)$ is a multiple, say $\alpha=p\tau$ for some $\tau\in H_4(G)$ and prime $p$. Then $$\chi(M)\ge 2 -2\  \dim_{\zz / p}(H^1(G;\zz / p)) + 2\ \dim_{\zz / p}(H^2(G;\zz / p)).$$

In particular, if $H_1(G;\zz)$ has no $p$--torsion, then $$\chi(M)\ge 2 -2\beta_1(G) + 2\beta_2(G).$$ \end{theorem}

\begin{proof}  Let $f\co M\to B_G$ denote the corresponding map.  The homology group $H_4(M;\zz / p)$ is isomorphic to $\zz / p$. The coefficient homomorphism $H_4(M;\zz)\to H_4(M;\zz / p)$ induced by the surjection $\zz\to \zz / p$ takes the fundamental class $[M]$ to a generator which we denote by $[M;\zz / p]$. Poincar\'e duality implies that the intersection form with $\zz / p$ coefficients $$H^2(M;\zz / p)\times H^2(M;\zz / p)\to \zz / p, \ (x,y)\mapsto (x\cup y)\cap [M;\zz / p]$$ is non-singular.

Since $\alpha=p\tau$, naturality of the coefficient homomorphism implies that $\alpha$ maps to zero under the morphism $f_*\co H_4(G;\zz)\to H_4(G;\zz / p)$. Thus the $\zz / p$ intersection form vanishes on the subspace $$I=\text{ Image}(f^*\co H^2(G;\zz / p)\to H^2(M;\zz / p)),$$ since $$0=(x\cup y)\cap 0=(x\cup y)\cap f_*([M;\zz / p])= ( f^*(x)\cup f^*(y))\cap [M;\zz / p].$$

The subspace $I$ has the same $\zz / p$--dimension  as $H^2(G;\zz / p)$; in fact we have that  $f^*\co H^2(G;\zz / p)\to H^2(M;\zz / p)$ is injective since $B_G$ can be constructed by adding cells of dimension 3 and higher to $M$.    Since the $\zz / p$--intersection form is non-singular, $$\dim_{\zz / p}H^2(M;\zz / p)\ge 2 \dim_{\zz / p} I.$$

Computing the Euler characteristic of $M$ using $\zz / p$ coefficients and $\zz / p$\,--Poincar\'e duality yields \begin{align*}\chi(M)&=  2 - 2 \dim_{\zz / p}H^1(M;\zz / p)+ \dim_{\zz / p}H^2(M;\zz / p)\\ &\ge 2 - 2 \dim_{\zz / p}H^1(G;\zz / p)+2 \dim_{\zz / p}H^2(G;\zz / p). \end{align*}

The last statement follows from the universal coefficient theorem. \end{proof}

\subsection{$L^2$ methods}

Several authors have explored applications of $L^2$ methods to the study of the geography problem of 4--manifolds.  Among these are Eckmann~\cite{e1,e2}, Kotschick~\cite{ko} and L\"uck~\cite{luck}.  An important result is  Theorem 5.1 
of~\cite{luck}, which states  that if $G$ is a group whose first $L^2$ Betti  number   vanishes, then any $M\in \calm(G)$ satisfies $\chi(M)\ge |\sigma(M)|$. Thus, for such a group $G$, $$ |\sigma|\leq q_G(\sigma).$$   Examples of such groups include amenable groups~\cite{cheeger-gromov}, extensions of finitely presented groups by $\zz$~\cite{jk, luck},  and fundamental groups of closed oriented 4--manifolds with a geometry (in the sense of Thurston) different from $S^2\times \hh^2$~\cite{luck}.

\bigskip


\section{Free Abelian Groups}

For the rest of this article we will focus on the geography of the free abelian groups. This class of groups provides a setting to investigate the ideas introduced above more deeply, and, as we shall see, will quickly lead us to difficult algebraic and 4--manifold problems.

\subsection{The homology and cohomology of $\zz^n$}

We begin by reviewing  the homology ring of $\zz^n$ and setting up notation and conventions.

A particular Eilenberg-MacLane space for $\zz^n$ is  the $n$--torus, $B_{\zz^n}=(S^1)^n=T^n$.  The cohomology ring $H^*(\zz^n)$ is  the exterior algebra on $H^1(\zz^n)$.   In particular, $H^k(\zz^n)$  is a free abelian group with rank the binomial coefficient, $C(n,k) $.   When convenient, we fix a basis $\gamma_1,\cdots, \gamma_n$ of $\pi_1(T^n)=H_1(\zz^n)$. We denote the dual basis by $x_1,\cdots, x_n\in H^1(T^n)$. We drop the notation ``$\cup$'' to indicate multiplication and so $H^k(\zz^n)$ has basis the $C(n,k)$ products $x_{i_1}x_{i_2}\cdots x_{i_k}, \ i_1<i_2<\cdots < i_k$.

The fact that $T^n$ is a closed orientable manifold gives us some additional algebraic structure which we now describe.  Note that the basis of $\pi_1(T^n)$ determines one of the two generators of $H_n(\zz^n)\cong\zz$ (that is, an orientation of $T^n$), and hence a fundamental class which we denote by $[T]\in H_n(\zz^n)$.  This then defines the duality isomorphism $$\cap [T]\co H^k(\zz^n)\to H_{n-k}(\zz^n).$$

Given $\alpha\in H_4(\zz^n)$ denote by $\omega\in H^{n-4}(\zz^n)$ the Poincar\'e dual to $\alpha$, so $$\omega\cap [T]=\alpha\in H_4(\zz^n).$$  Given $x,y\in H^2(\zz^n)$, $$(xy)\cap \alpha= (xy\omega)\cap [T].$$ For our purposes, this formula is best recast in the following proposition, whose proof is just an application of the formula of Equation (\ref{eq5}).

 \begin{prop}\label{propzn} Let $\alpha\in H_4(\zz^n)$ and choose $(M,f)\in \calm(\zz^n,\alpha)$. Let $\omega\in H^{n-4}(\zz^n)$ satisfy $\omega\cap[T]=\alpha$.

Then the restriction of the intersection form $H^2(M)\times H^2(M)\to \zz$  of $M$ to the subgroup $I(M,f)=\text{\rm  image} (f^*\co H^2(\zz^n)\to H^2(M))$  is given by the pairing $$H^2(\zz^n)\times H^2(\zz^n)\to \zz, \ (x,y)\mapsto xy\omega\cap [T].$$ 
\end{prop} 
 \qed
 
The pairing of Proposition \ref{propzn} does not depend on the orientation of $T^n$, since the intersection form on $M$ is independent of the orientation of $T^n$. In any case it is useful to observe that  changing the orientation of $T^n$ (for example, by changing the basis $\gamma_1,\cdots, \gamma_n $ of $\pi_1(M)$) changes the signs of both $\omega$ and $[T]$.

Notice the simplicity of this pairing: $H^*(\zz^n)= \Lambda^*(\zz^n)$, so an orientation of $T^n$ is just an identification $\Lambda^n(\zz^n)\cong \zz$, and the pairing of Proposition \ref{propzn} in the usual notation of the exterior algebra is \begin{equation}\label{ext} \Lambda^2(\zz^n)\times \Lambda^2(\zz^n)\to \zz, \ (x,y)\mapsto x\wedge y\wedge \omega.\end{equation}

Some basic consequences of this observation are assembled in the following theorem.

\begin{theorem}\label{exterior}  Fix $\alpha\in H_4(\zz^n)$ and $(M,f)\in \calm(\zz^n,\alpha)$. Use the  isomorphism $f^*\co H^1(\zz^n)\to H^1(M)$ to identify the basis $\{x_i\}$ of $H^1(\zz^n)$ with a basis of $H^1(M)$.

Then the subspace $I(M,f)\subset H^2(M)$ is a free abelian summand of rank $C(n,2)$.  In the basis $\{x_ix_j\}_{i<j}$ of $I(M,f)$ one has $(x_i x_j )(x_k x_\ell)=0$ if an index is repeated in the set $\{i,j,k,\ell\}$. In particular
 \begin{enumerate} 
\item The restriction of the intersection form of $M$ to $I(M,f)$ is even. 
\item The intersection form of $M$ contains  $n-1$ dimensional isotropic subspaces,  for instance  the subspace spanned by $\{x_1x_2, x_1x_3,\cdots, x_1x_n\}$. Hence $q_{\zz^n}(\sigma)\ge |\sigma|$. 

\item $(x_ix_j)(x_kx_\ell)=-(x_ix_k)(x_jx_\ell)$. \end{enumerate}

\end{theorem} \begin{proof}  Since $f^*\co H^2(\zz^n)\to H^2(M)$ is injective by Theorem \ref{thm1},  $I(M,f)\subset H^2(M)$ is  free abelian of rank $C(n,2)$. To see that it is a summand, consider a map on the $n$--fold wedge $$j \co S^1\vee\cdots \vee S^1\to M$$ taking the $i$th circle to a loop representing the generator $\gamma_i\in \pi_1(M)$. Since $\pi_1(M)=  \zz^n$ is abelian, the map $j$ extends to the 2--skeleton of the $n$--torus $$j \co (T^n)^{(2)}\to M.$$ The induced composite on cohomology $$H^2((T^n)^{(2)})\cong H^2(T^n)\xrightarrow{f^*}H^2(M) \xrightarrow{j^*} H^2((T^n)^{(2)})$$ is clearly the identity and hence gives a splitting of $f^*$.

The rest of proof follows from Corollary  \ref{cor1}  and basic properties of the exterior algebra, such as $x_i x_j=-x_jx_i$. In particular, since $(x_ix_j)(x_ix_j)=0$, the intersection form of $M$ restricted to $I(M,f)$ has  zeros on the diagonal in this basis and hence is even.  Notice  that since the intersection form of $M$ has an $(n-1)$ dimensional isotropic subspace, $\beta_2(M)\ge |\sigma| + 2(n-1)$, which implies that $\chi(M)\ge |\sigma|$.  \end{proof}

\begin{corollary} \label{cor6.4}Suppose that $\alpha\in H_4(\zz^n)$ and that there exists an $(M,f)\in \calm(\zz^n,\alpha)$  so that $f^*\co H^2(\zz^n)\to H^2(M)$ is surjective.
Then $n \ne 2,3, $ and the intersection form of $M$ is equivalent to $k E_8 \oplus \ell H$ for some $k\in \zz$ and $\ell \ge n-1$. In particular $C(n,2)=8 |k|  + 2\ell$ and therefore is even. If $M$ is smooth, then in addition $|k|$ is even. \end{corollary}

\begin{proof} If $f^*\co H^2(T^n)\to H^*(M)$ is surjective, then by Theorem \ref{exterior} it is an isomorphism and the intersection form of $M$ is even. Since it contains an $n-1$ dimensional isotropic subspace, $C(n,2)\ge 2(n-1)$, and hence $n \ne 2,3$. Furthermore the intersection form is indefinite, so the classification of even, unimodular integer bilinear forms shows that the intersection form of $M$ is equivalent to   $k E_8 \oplus \ell H$ for some $k\in \zz$ and $\ell \ge n-1$.

If $M$ is smooth, then since $H_1(M)$ has no 2--torsion it follows that $M$ is spin, so that by Rohlin's theorem the signature of $M$ is divisible by 16. Hence $k$ is even. \end{proof}

Theorem \ref{exterior} and  Corollary  \ref{cor6.4} together imply that $$q(\zz^n)\ge 2-2n + C(n,2) +\epsilon_n,$$ where $\epsilon_n$ equals zero if $C(n,2)$ is even and  equals $1$ if $C(n,2)$ is odd. The main result of ~\cite{kl} is that this lower bound is achieved for all $n$, except that $q(\zz^3) = 2$, rather than 0, and $q(\zz^5) = 6$, rather than $2$.  Moreover, the  examples constructed in \cite{kl} which realize this lower bound all have signature zero, since they are obtained by surgeries on connected sums of products of orientable surfaces. Thus we conclude that for $n\ne 3,5$, \begin{equation}\label{eq5.1}q_{\zz^n}(\sigma)\ge q_{\zz^n}(0)=2-2n + C(n,2) +\epsilon_n\end{equation}

\medskip

\noindent{\bf Example.} Suppose $\omega\in H^4(\zz^8)$ so that the induced pairing (\ref{ext}) is unimodular. Notice that  such an $\omega$ exists by the results of~\cite{kl} since $C(8,2)=28$ is even. The classification theorem for unimodular integer forms~\cite{MH} shows that this pairing is equivalent (after perhaps changing orientation) to either $14H$ or  $E_8\oplus 10H$.

Now if $(M,f)\in\calm(\zz^8, \omega\cap[T])$ satisfies $\beta_2(M)=28$, then $I(M,f)=H^2(M)$ and so intersection form of $M$ is equivalent to (\ref{ext}). The example constructed 
in~\cite{kl} has signature zero, hence has intersection form $14H$, but Corollary \ref{cor6.4} allows the possibility that for some $\omega$ the signature equals 8.  

\medskip \noindent{\bf Question.}  Does there exist $\omega\in H^4(\zz^8)$ so that (\ref{ext}) is equivalent to $E_8\oplus 10H$?  If so, does there exist $(M,f)\in \calm(\zz^8,\omega\cap[T])$ so that $\beta_2(M)=28$ and $f_*([M])=\omega\cap[T]$?  Note that by Rohlin's theorem  such an $M$ cannot be smooth.

\subsection{Calculations}

\subsubsection {\bf  n = 0,1,2,3 }   Since $\zz^0,\zz^1$ and $ \zz^3$ are   closed 3--manifold groups (the corresponding 3--manifolds are $S^3$, $S^1\times S^2$, and $T^3$) and $\zz^2$ is the fundamental group of a torus,  the geography for these groups is given by Theorems \ref{3manifolds} and \ref{surface}: \[ q_{\zz^n}(\sigma)= \begin{cases} |\sigma| + 2 & \text{if } n=0,3\\ |\sigma|  & \text{if } n=1, 2 \end{cases}\]

\begin{figure}[h] \centerline{ \fig{.4}{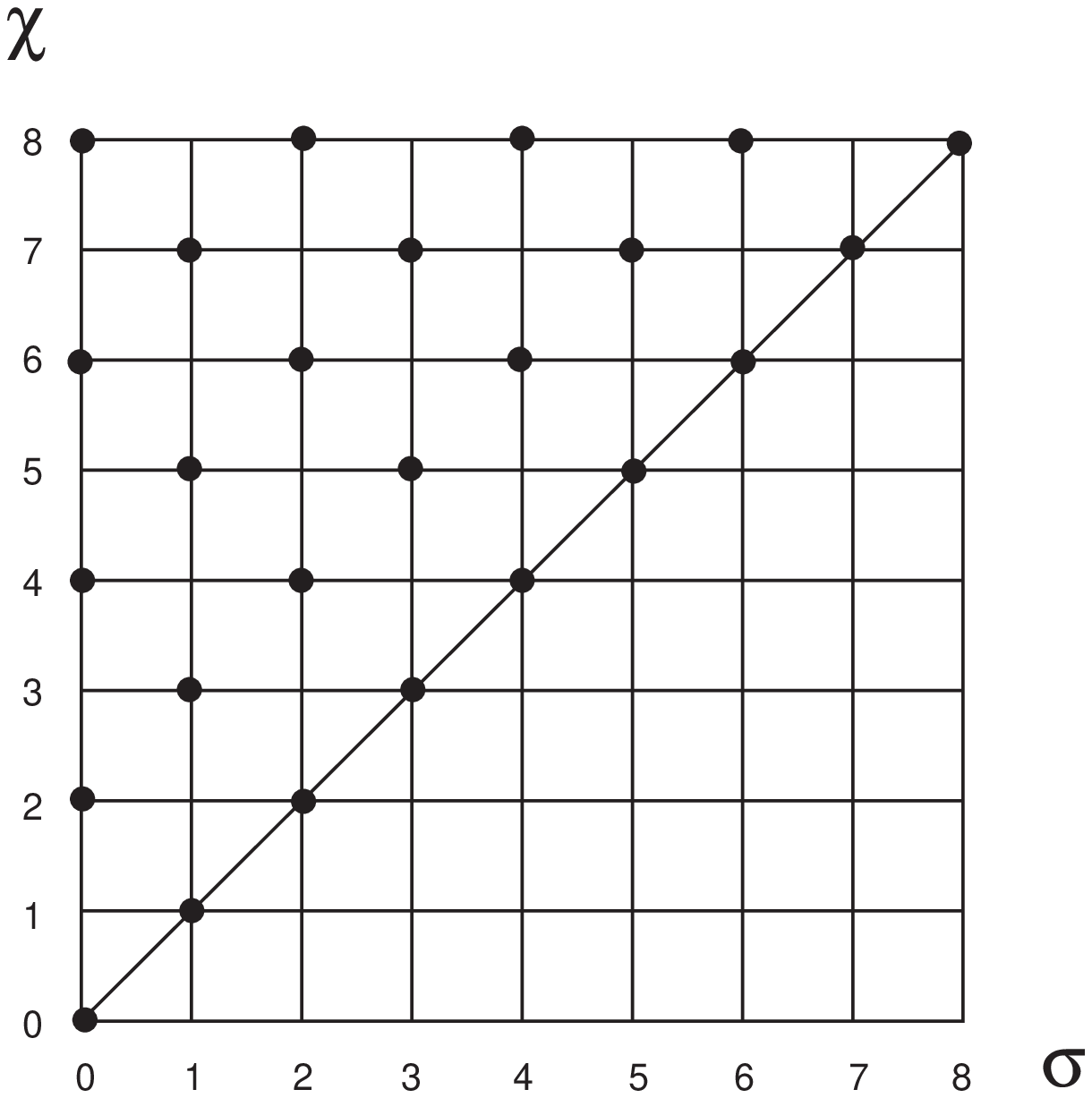}} \caption{The geography of $\zz$, $\zz^2,$ and $\zz^4$}\label{q=s} \end{figure}

The geography of $\zz^0$ and $\zz^3$ has already been illustrated in Figure~\ref{q=0} and the geography of $\zz^1, \zz^2$,  and (as we will show next), $\zz^4$ is illustrated in Figure~\ref{q=s}.

\subsubsection{\bf n = 4.}  First note that $H_4(\zz^4)=\zz$. Moreover, since the 4--torus admits an orientation reversing homeomorphism, $q_{\zz^4,-\alpha}(\sigma)=q_{\zz^4,\alpha}(\sigma)$ for any $\alpha\in H_4(\zz^4)$.

Since $\beta_2(\zz^4)=C(4,2)=6$, Theorem \ref{exterior}  implies that for any $\alpha\in H_4(\zz^4)$ and any $(M,f)\in \calm(\zz^4,\alpha)$, $\beta_2(M)\ge 6$.  Theorem \ref{exterior}  shows that $q_{\zz^4,\alpha}(\sigma)\ge |\sigma|$. Taking the identity map of the 4--torus and connected sums with $\pm \cc P^2$ shows that $$q_{\zz^4}(\sigma)=|\sigma|$$ and that if $[T]\in H_4(\zz^4)$ is a generator $$q_{\zz^4,[T]}(\sigma)=|\sigma|.$$

If $\alpha=k[T]$ for $|k|>1$ or $k=0$, then Theorem \ref{multiple} implies that $\chi(M)\ge 2-8 + 12=6$ and so $q_{\zz^4,\alpha}(\sigma)\ge \max\{6,|\sigma|\}.$

A signature zero example with $\chi=6$ representing $\alpha=k[T]$, $|k|\ne 1$, can be constructed as follows.  Let $M':=T^4\# (S^1\times S^3)$, with $\pi_1(T^4)$ generated by $a_1,a_2,a_3,a_4$ and $\pi_1(S^1\times S^3)$ generated by $b$. Let $f'\co M'\to T^4$ be a map inducing the map $a_1\mapsto\gamma_1, a_2\mapsto\gamma_2,a_3\mapsto\gamma_3,a_4\mapsto\gamma_4^k, b\mapsto\gamma_4$ on fundamental groups.  Then $f'_*([M'])=k[T]$ and   $f'$ induces an epimorphism on fundamental groups. Note that $\chi(M')=-2$ and $\sigma(M')=0$. Now perform 4 surgeries along circles in $M'$: the first  to introduce the relation $b^k=a_4$ and the other three to introduce the commutator relations $[a_1,b],[a_2,b]$ and $[a_3,b]$. Each surgery increases the Euler characteristic by 2 and leaves the signature invariant. The map $fÕ$ clearly extends over the trace of the surgery. Thus the resulting 4--manifold $M$ has a map $f \co M\to T^4$ inducing an isomorphism on fundamental groups, and representing $k[T]$. Moreover, $\chi(M)=6$ and $\sigma(M)=0$.

Thus $ q_{\zz^4,k[T]}(0)=6$ for $|k|\ne 1$. Taking connected sums with $\pm \cc P^2$ shows that $$ \max\{6,|\sigma|\}\leq q_{\zz^4,k[T]}(\sigma)\leq |\sigma| +6.$$

When  $k=0$, i.e.~$(M,f)\in\calm(\zz^4,0)$, then $I(M,f)$ is a 6--dimensional isotropic subspace. Hence \begin{equation*} q_{\zz^4,0}(\sigma)= |\sigma| +6. \end{equation*}

We summarize these calculations in the following theorem. 

\begin{theorem}\label{nis4}  Given $\alpha\in H_4(\zz^4)$ define the non-negative integer  $k\ge 0$ so that $\pm \alpha=k[T]$. Then the following hold. \begin{enumerate} \item $q_{\zz^4}(\sigma)=|\sigma|$, and so $q(\zz^4)=0$, $p(\zz^4)=0$. \item If $k=0$,  $q_{\zz^4,\alpha}(\sigma)=|\sigma|+ 6$. \item If $k=1$,   $q_{\zz^4, \alpha}(\sigma)=|\sigma|.$ \item If $k>1$,  $ \max\{6,|\sigma|\}\leq q_{\zz^4,\alpha}(\sigma)\leq |\sigma| +6.$ \end{enumerate}\qed \end{theorem}

\noindent{\bf Question.} Is $q_{\zz^4,k[T]}(\sigma)=|\sigma|+6$ for $|k|>1$? From parity considerations we know that $q_{\zz^4,k[T]}(1)=7$ when $|k|\ne 1$. We do not know whether $q_{\zz^4,k[T]}(2)=6$ or  $8$.

\begin{figure}[h]
 \centerline{ \fig{.4}{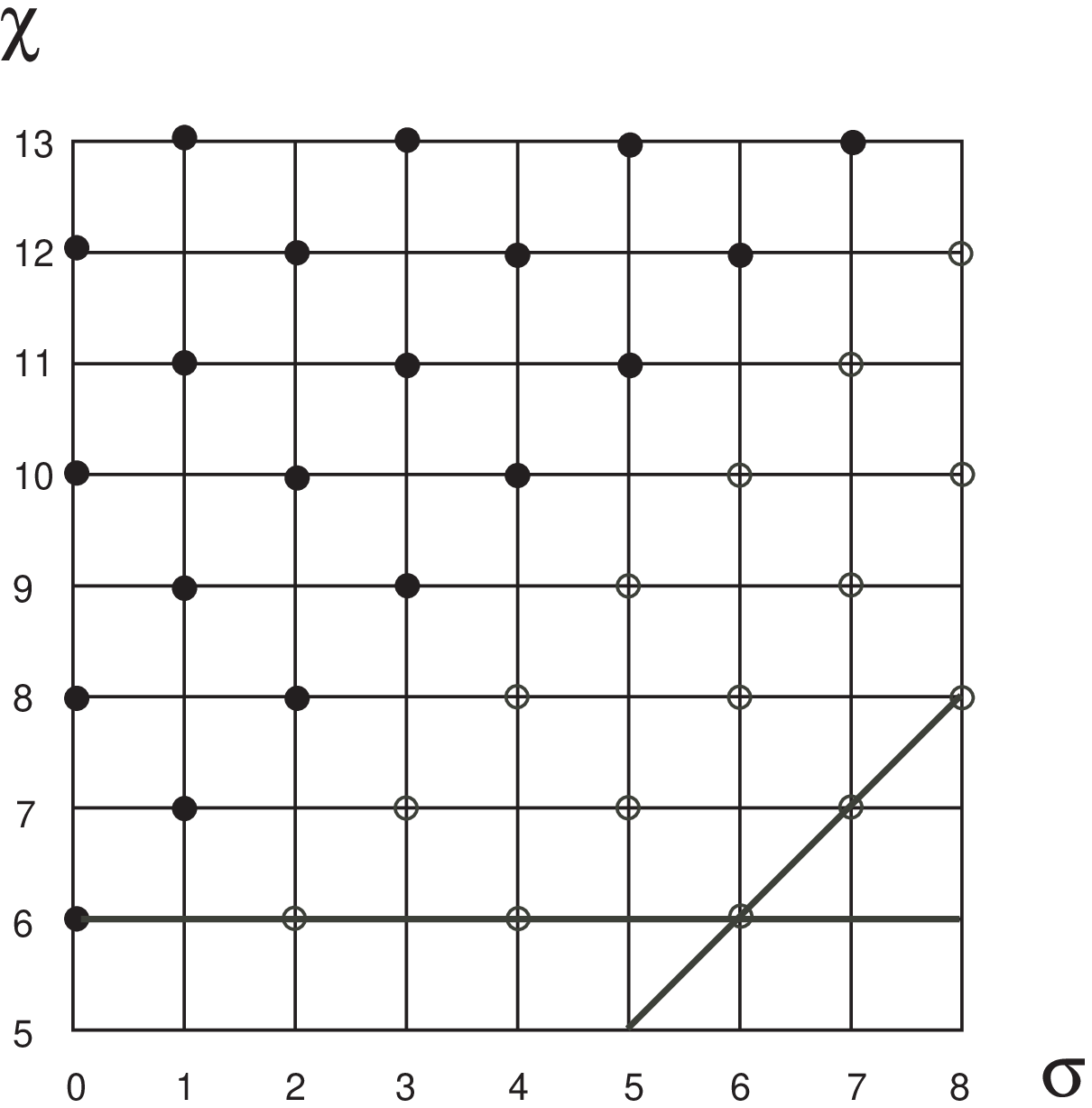}} 
\caption{The unknown points 
in $\calg(\zz^4, k[T]) , |k| \ge 2$} \label{unknown4} 
\end{figure}

The unknown values of $q_{\zz^4, k[T]}$ for $|k| >1 $ are marked with circles in 
Figure~\ref{unknown4}.

\subsubsection{\bf n = 5.}   Note that $H^2(\zz^5)$ has rank 10 and $H_4(\zz^5)$ has rank $5$. We first claim that given any $\alpha\in H_4(\zz^5)$, there is an automorphism $\Psi\in\Aut(\zz^5)$ so that $\Psi_*(\alpha)= (kx_1)\cap [T]$ for some $k\in \zz$.

The following lemma helps to keep track of how automorphisms of $\pi_1(T^n)=H_1(T^n)$ act on the cohomology $H^*(T^n)$, thinking of the cohomology as the exterior algebra.

\begin{lemma}\label{changebasis}  Let $A \co \pi_1(T^n)\to \pi_1(T^n)$ be an automorphism, given in the basis $\gamma_1\cdots, \gamma_n$ as an $n\times n $ matrix $A\in GL(\zz^n)$. Let $x_1,\cdots,x_n\in H^1(\zz^n)=\Hom(\pi_1(T^n),\zz)$ denote the dual basis. Then the induced map $A^* \co H^1(\zz^n)\to H^1(\zz^n)$ is given by the transpose $A^T$ in the basis $\{x_i\}$.

Moreover,   if $\omega\in H^{n-4}(\zz^n)$, $\alpha\in H_4(\zz^n)$ satisfy $\omega\cap [T]=\alpha$, then $A_*(\alpha)= \det(A)((A^*)^{-1}(\omega))\cap [T]$. \end{lemma}

\begin{proof} The fact that $A^*=A^T$ is well-known and easy. The other formula is a consequence of the the naturality of cap products, expressed by the identity $f_*( f^*(x)\cap y)=x \cap f_*(y)$, valid for any continuous map $f$.  In more detail, taking $f:T^n\to T^n$ to be a map inducing the automorphism $A:\pi_1(T^n)\to \pi_1(T^n)$,    the induced map $A_* \co H_n(\zz^n)\to H_n(\zz^n)$ is multiplication by $\det(A)$. One computes 
 \begin{align*} A_*(\alpha)&= A_*(\omega\cap [T])=A_*(A^*(A^*)^{-1}(\omega)\cap [T])\\ &= (A^*)^{-1}(\omega)\cap A_*([T])= (A^*)^{-1}(\omega)\cap \det(A)[T]. \end{align*} \end{proof}

\medskip

We use Lemma \ref{changebasis} as follows. Given $\alpha\in H_4(\zz^5)$, let $\omega\in H^1(\zz^5)$ satisfy $\omega\cap [T]=\alpha$. Since $\omega\in H^1(\zz^5)$, we can write $\omega= k(b_1x_1+\cdots + b_5x_5)$ where $k\in\zz$ and $b_1x_1+\cdots + b_5x_5$ is primitive (that is, not divisible). Thus there is some matrix $B\in GL(\zz^5)$ which sends $\omega$ to $k x_1$. Clearly we can assume $\det(B)=1$ and $k\ge 0$. Lemma \ref{changebasis} implies that the matrix $A= (B^T)^{-1}$ determines an automorphism $\Psi\co  \zz^5 \to  \zz^5$ satisfying  $\Psi_*(\alpha)= kx_1\cap [T]$.

Since $\calg(\zz^5,\alpha)=\calg(\zz^5,\Psi_*(\alpha))$, we see that  $\calg(\zz^5,\alpha)=\calg(\zz^5, kx_1\cap [T])$, where $k\ge 0$ is the unique integer so that $\alpha=k\tau$ for $\tau$ a primitive class.

If $\alpha=0$ then $I(M,f)$ is a 10-dimensional isotropic subspace for  any $(M,f)\in \calm(\zz^n,0)$.  For such $(M,f)$,  $|\sigma|\leq \beta_2(M)-20$, and so $$q_{\zz^5,0}(\sigma)=2-10+\beta_2(M)\ge |\sigma|+ 12.$$

If $\alpha\ne 0$, then since $xy\cap\alpha=xy\omega\cap [T^5]=kxyx_1\cap[T]$, the intersection form vanishes on the subspace of  $I(M,f))$ generated by 

\begin{equation}\label{5iso}
\{x_1x_2, x_1x_3, x_1x_4, x_1x_5, x_2x_3, x_2x_4, x_2x_5\},
\end{equation}
 which is of rank 7.  Thus,  for any pair  $(M,f) \in \calm(\zz^5,\alpha)$ with $\alpha\ne 0$,  $\chi(M)   =  2  -10 + \beta_2(M) \ge -8 + |\sigma(M)| + 14 = |\sigma(M)| +6$. Hence for $\alpha\ne 0$, $$ q_{\zz^5,\alpha}(\sigma)\ge |\sigma|+6.$$

To construct an upper bound, take $M'=T^4\#(T^2\times S^2)$ with $\pi_1(T^4)$ generated by $a_1,a_2,a_3, a_4 $ and $\pi_1(T^2\times S^2)$ generated by $b_1,b_2$. Thus $\chi(M')=-2$ and $\sigma(M')=0$. Select a map $M'\to T^5$ inducing the homomorphism $a_1\mapsto \gamma_2, a_2\mapsto \gamma_3, a_3\mapsto \gamma_4, a_4\mapsto \gamma_5^k, b_1\mapsto \gamma_1, b_2\mapsto \gamma_5$ on fundamental groups. Then $\alpha=kx_1\cap [T]$.

If $\alpha=0$ (in other words $k=0$), perform 7 surgeries on $M'$: one to kill $a_4$  and six to kill the commutators of $a_1,a_2,a_3$ with $b_1,b_2$. Each surgery increases $\chi$ by $2$. Clearly the map from the resulting manifold $M$ to $T^5$ induces an isomorphism and represents $0$ in $H_4(\zz^5)$. Since $\chi(M)=12$ and $\sigma(M)=0$, taking connected sums with $\pm \cc P^2$ provides examples showing $q_{\zz^5,0}(\sigma)\leq |\sigma|+ 12,$ and so $$q_{\zz^5,0}(\sigma)= |\sigma|+ 12.$$

If $\alpha=x_1\cap [T]$, i.e.~$k=1$, perform 4  surgeries on $M'$: one to kill $b_2a_4^{-1}$, and three to kill the commutators $[a_1,b_1], [a_2, b_1]$ and $[a_3,b_1]$. The map from the resulting manifold $M$ to $T^5$ induces an isomorphism and represents $\alpha $ in $H_4(\zz^5)$. Since $\chi(M)=6$ and $\sigma(M)=0$, taking connected sums with $\pm \cc P^2$ provides examples showing $q_{\zz^5,0}(\sigma)\leq |\sigma|+ 6$. Thus for $k=1$: $$q_{\zz^5,x_1\cap[T]}(\sigma)= |\sigma|+ 6.$$

If $k>1$, perform 7  surgeries on $M'$: one to kill $b_1^ka_4^{-1}$, and  six  to kill the commutators $[a_1,b_1], [a_2,b_1],[a_3,b_1],   [a_1,b_2], [a_2,b_2]$, and $[a_3, b_2]$. The map from the resulting manifold $M$ to $T^5$ induces an isomorphism and represents $\alpha $ in $H_4(\zz^5)$. Since $\chi(M)=12$ and $\sigma(M)=0$, taking connected sums with $\pm \cc P^2$ provides examples showing $q_{\zz^5,kx_1\cap [T]}(\sigma)\leq |\sigma|+ 12$.   On the other hand, Theorem \ref{multiple} shows that for such $\alpha$, $\chi(M)\ge 12$.

Thus for $k>1$: $$\max\{12, |\sigma|+ 6 \} \leq q_{\zz^5,kx_1\cap[T]}(\sigma)\le  |\sigma|+ 12.$$

Taking the minimum over all homology classes $\alpha$ we see that $$q_{\zz^5}(\sigma)=|\sigma| + 6.$$

In summary: \begin{theorem}\label{nis5}  Given $\alpha\in H_4(\zz^5)$, there exists a unique non-negative integer $k$  and $A\in \Aut(\zz^5)$ so that $A_*(\alpha)=(kx_1)\cap [T]$. Then the following hold. \begin{enumerate} \item $q_{\zz^5}(\sigma)=|\sigma|+6$, and so $q(\zz^5)=6$, $p(\zz^5)=6$. 
\item If $k=0$, $q_{\zz^5,0}(\sigma)=|\sigma|+ 12$. \item If $k=1$  $q_{\zz^5, x_1\cap [T]}(\sigma)=|\sigma|+6.$ 
\item If $k>1$, $ \max\{12,|\sigma|+6\}\leq q_{\zz^5,kx_1\cap[T]}(\sigma)\leq |\sigma| +12.$ \end{enumerate}\qed \end{theorem}

\noindent{\bf Question.} If $k>1$, is $q_{\zz^5,kx_1 \cap [T] }(\sigma)= |\sigma|+ 12$?

\medskip

\begin{figure}[h] \centerline{ \fig{.4}{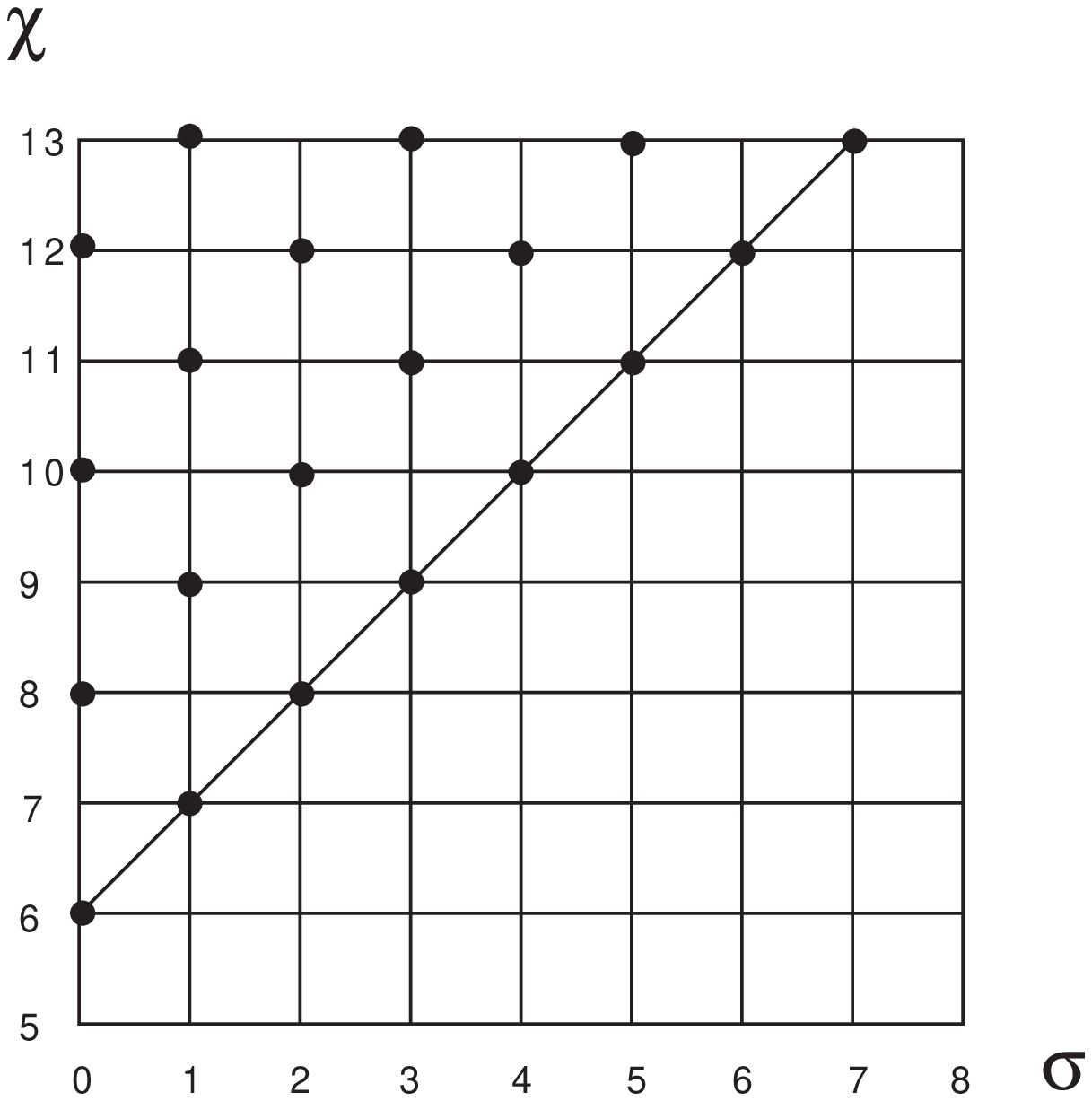}} \caption{The geography of $\zz^5$}\label{q=5} \end{figure}

The geography of $\zz^5$ is illustrated in Figure~\ref{q=5}. (Note the change in the vertical scale in the figure.)

\subsubsection{\bf n = 6.}\label{sub6}     For $\zz^6$, the situation becomes more delicate and interesting; for instance, $q_{\zz^6}(\sigma)$ is not linear for $\sigma \ge 0$.   Notice first that since $\zz^6$ maps onto $\zz^5$, the metabolizer found in the $n=5$ case above yields a metabolizer of rank $7$ for any 4--manifold $M$ with $\pi_1(M) \cong \zz^6$.  Thus,  if  $\pi_1(M)=\zz^6$ then $$\chi(M) = 2 - 12 + \beta_2(M) \ge -10 + |\sigma(M)| + 14 = 4 + |\sigma(M)|.$$  Moreover, since $C(6,2) = 15$ is odd, we also have that $\beta_2(M) \ge 16$ by Corollary~\ref{cor6.4}, so $\chi(M) \ge 6$, and hence $$\max\{6,|\sigma|+4\}\leq q_{\zz^6}(\sigma).$$

It is convenient at this point to recall the examples of symmetric products,~\cite{{Baldridge-Kirk}, Mc}. 

\begin{prop}\label{sym} Let $F_k$ denote the oriented compact surface of genus $k$. Let $S_{2k}=Sym^2(F_k)$. Then $S_{2k}$ is a smooth 4--manifold with $\pi_1(S_{2k})=\zz^{2k}$, $\chi(S_{2k})=2k^2-5k+3$ and $\sigma(S_{2k})=1-k$. \qed \end{prop}

\begin{figure}[h] \centerline{ \fig{.4}{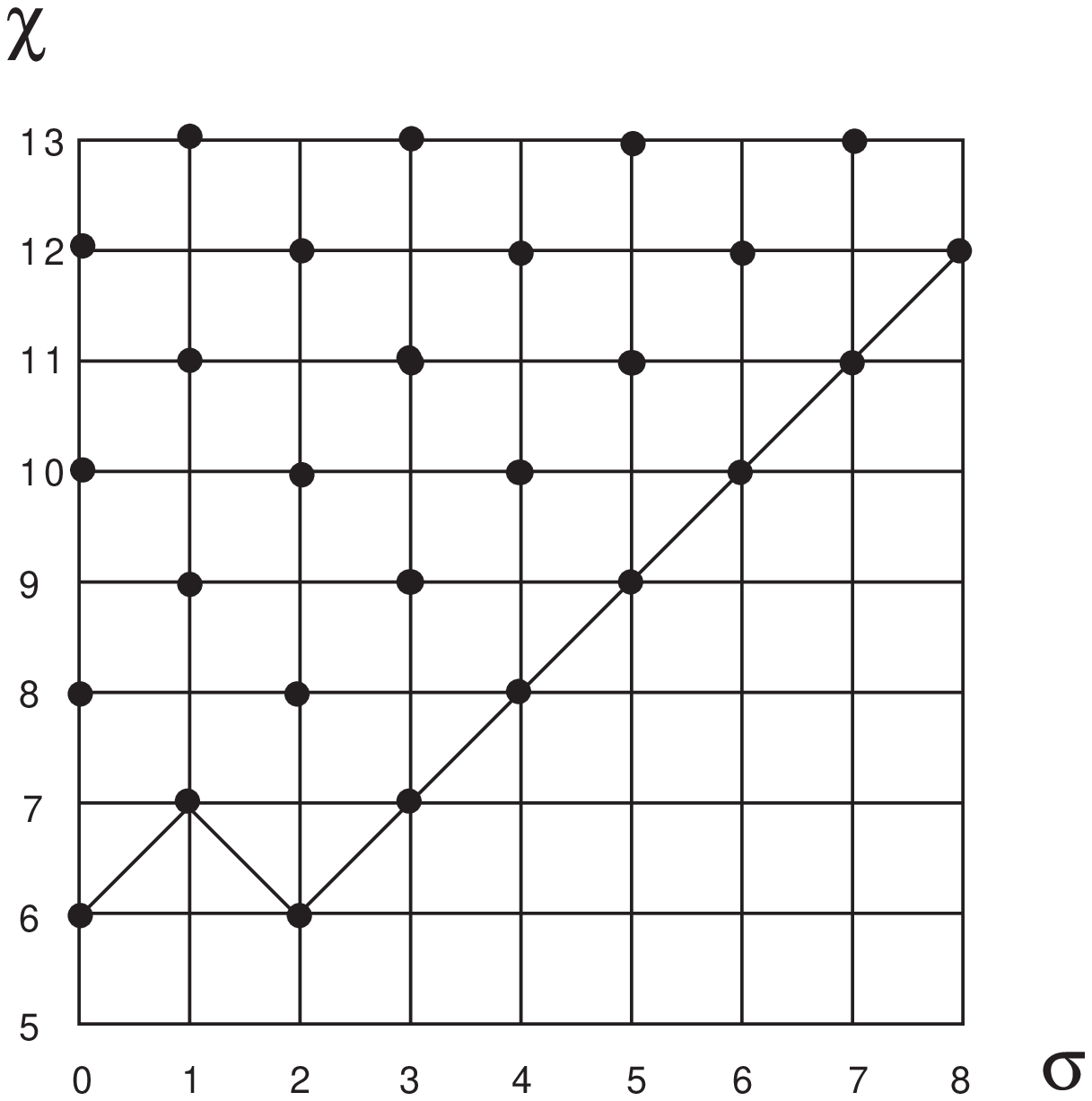}} \caption{The geography of $\zz^6$}\label{q=6} \end{figure}

A 4--manifold $M$ is constructed in~\cite{kl}    with $\pi_1(M) = \zz^6$,  $\chi(M) = 6$ and $\sigma(M) = 0$ (see also the proof of Theorem \ref{nequals6}  below).    It follows that $q_{\zz^6}(0) = 6.$  Since   $q_{\zz^6}(1) \ge 6$, it follows that  $q_{\zz^6}(1) = 7.$   The 4--manifold $S_{6}$ has fundamental group $\zz^6$, $\chi(S_6)=6$, and $\sigma(S_6)=-2$.   Thus, $q_{\zz^6}(2)=q_{\zz^6}(-2) = 6.$  It follows that \[ q_{\zz^6}(\sigma)=\begin{cases} 6 & \sigma=0,\\ 7 & |\sigma|=1,\\ |\sigma|+4 &|\sigma|\ge 2,\end{cases}\] as illustrated in Figure~\ref{q=6}.

To compute $q_{\zz^6,\alpha}$,  we find a canonical form for   homology classes $\alpha\in H^2(\zz^6)$.   In the following theorem, we fix a basis $\gamma_1,\dots, \gamma_6\in H_1(\zz^6) $ and let $x_1,\cdots,x_6\in H^1(\zz^6)$ denote the dual basis. These determine an orientation class $[T]\in H_6(\zz^6)$.

For a triple of integers $a,b,c$, the notation $a|b|c$ means that there are integers $k,\ell$ so that $b=ak$ and $c=b\ell$.

\begin{theorem} \label{6normalform} Given any  $\alpha\in H_4(\zz^6)$, there exists an automorphism  $A\in \Aut(\zz^6)$ and a unique set of integers $a,b,c\in \zz$ with $a|b|c$,  $a, b \ge 0$, so that $$A_*(\alpha)= (a  x_1x_2 + bx_3x_4 + cx_5x_6)\cap [T].$$ 
\end{theorem}

 \begin{proof} Let $B\in SL_6(\zz)$.   In the basis $\{x_i\}$,  $B$ defines automorphisms $B^*\co H^1(\zz^6)\to H^1(\zz^6)$ and $B^*\co H^2(\zz^6)\to H^2(\zz^6)$ by $B^*(x_i)= \sum_j B_{ij}x_j$ and   $B^*(x_ix_j)=B^*(x_i)B^*(x_j)$.
Let $\omega\in H^2(\zz^6)$ satisfy $\omega\cap[T]=\alpha$. Lemma \ref{changebasis} implies that if we set $A=(B^T)^{-1}$, the automorphism $A\co \zz^6\to \zz^6$ given by $A(\gamma_i)=\sum_j A_{i,j}\gamma_j$ satisfies $A_*(\alpha)=B^*(\omega)\cap[T]$. 

We first prove that for any  $\omega\in H^2(\zz^6)$, there is a matrix $B\in SL_6(\zz)$ so that $B^*(\omega)$ is of the form $a  x_1x_2 +  bx_3x_4 +  cx_5x_6$ for some $a, b, $ and $c$.

Any  $\omega \in H^2(\zz^6)$ can be written as $\sum \alpha_{i,j} x_i x_j,  i <j$. Collecting terms, we can express $\omega$ as $$\omega = a_1x_1(\alpha_2 x_2 + \cdots +  \alpha_6 x_6) + \omega_2,$$ where the $\alpha_i$ have $\gcd =1$, $a_1\ge 0$,  and $\omega_2$ can be expressed in terms of $x_i x_j$ with $2 \le i < j$.   (If $x_1$ appears in $\omega$, then   $a_1 > 0$ and the $\alpha_i$ are uniquely determined; otherwise, we have the  decomposition $\omega = 0x_1(x_2) + \omega_2$.)

Within $\langle x_2, \ldots , x_6\rangle  \cong \zz^5$ there is a basis with its first element $\alpha_2 x_2 + \cdots +\alpha_6 x_6$.  Renaming those basis elements $x_2, \ldots , x_6$ gives $$ \omega = a_1x_1x_2 + \omega_2,$$ where $\omega_2$ continues to be expressed in terms of $x_ix_j, 2 \le i <j$.

Repeating this process,  we can change basis so that $$\omega = a_1 x_1x_2 + a_2 x_2x_3 + a_3x_3x_4 +a_4x_4x_5 + a_5 x_5 x_6 $$ and by changing signs of the appropriate $x_i$ we can arrange that $a_i\ge 0 $ for all $i$.

Now we want to see that by a further change of basis it can be assumed that either $a_1 = 0$ or $ a_2 = 0$.  We will demonstrate this by repeatedly changing basis so that $\min \{a_1, a_2\}$ is reduced.

Suppose now that $a_1 \le a_2$.  Make the change of basis $x_1 \mapsto x_1 + x_3$.    Then we have: $$\omega =  a_1x_1x_2  + (a_2 - a_1)x_2x_3 + a_3x_3x_4 +a_4x_4x_5 + a_5 x_5 x_6.$$ The coefficient on $x_2x_3$ became smaller, and this can be repeated until the coefficient on $x_2x_3$ is smaller than $a_1$.

On the other hand, if $a_2 < a_1$, make the change of basis $x_2 \mapsto x_1 + x_2$.  Then we have $$\omega =  x_1(a_1x_2  + a_2x_3)  + a_2x_2x_3 + a_3x_3x_4 +a_4x_4x_5 + a_5 x_5 x_6.$$ Now, factor out $a_1' = \gcd(a_1,a_2)$ (and note that $a_1' \le a_2$) to write $$\omega =  a_1'x_1(\alpha_1x_2  + \alpha_2x_3)  + a_2x_2x_3 + a_3x_3x_4 +a_4x_4x_5 + a_5 x_5 x_6,$$ with $\alpha_1$ and $\alpha_2$ relatively prime.  As earlier, a change of basis on $\langle x_2, \ldots , x_6\rangle$   can be made so that $\alpha_2 x_2 + \alpha_3x_3$ is replaced by $x_2$ and with respect to this new basis, we have $$\omega =  a_1'x_1x_2  + a_2'x_2x_3 + a_3'x_3x_4 +a_4'x_4x_5 + a_5' x_5 x_6.$$ Notice that $a_2' $ may be greater than $a_2$, but we do have that $a_1' \le a_2 < a_1$. If $a_1' < a_2$ we are done.  Otherwise, do the earlier change of basis to  lower $a_2'$ by a multiple of $a_1'$ so that $a_2' < a_1'$ and thus  $a_2' < a_1' \le a_2 < a_1$.  Hence the smaller of the two decreased.   Continuing in this manner we obtain either $a_1=0$ or $a_2=0$.

With this, it follows from induction that $\omega \in H^2(\zz^6)$ can be put in the form $\omega = ax_1x_2 + bx_3x_4 + cx_5x_6$ for some $a, b$, and $c$, with respect  to some basis.  
 
 We next want to achieve the desired divisibility of the coefficients.
Consider $c(ax_1x_2 + b x_3x_4)\in H^2(\zz^4)$ with $a,b,c\in \zz$ and  $a$ and $b$ relatively prime. Choose integers $p,q$ with $ap+bq=1$. The determinant 1 automorphism defined by the substitutions $$x_1= x_1'-bqx_3', x_2=px_2'-b x_4', x_3=x_1'+ap x_3', x_4=qx_2'+ ax_4'$$ takes $c(ax_1x_2 + b x_3x_4)$ to $c(x_1'x_2'+ ab x_3'x_4')$.
With this, it follows from induction that $\omega \in H^2(\zz^6)$ can be put in the form $\omega = ax_1x_2 + bx_3x_4 + cx_5x_6$ with respect  to some basis, with $a|b|c$.  

Changing the sign of basis elements can then ensure that  $a,b,c\ge 0$. We can further arrange  that  the basis change  matrix has determinant $1$,  perhaps at the cost of  changing the sign  of $c$.  The remark at the start of the proof then  implies that there is an automorphism $A$ with $A_*(\alpha)= (ax_1x_2 + bx_3x_4 + cx_5x_6)\cap [T]$, with $a,b\ge 0$ and $a|b|c$.

To show the uniqueness of $a, b$, and $c$ we consider the quotient  of $H^2(\zz^6) / \langle \omega \rangle $, now viewed as an  abelian group.  The order of the torsion of the quotient group is $a$.  Thus, the value of $a$ depends only on the $ \Aut(\zz^6)$ orbit of $\omega$.

Next   consider  $\omega^2\in H^4(\zz^6)$. When placed in normal form, $$\omega^2= 2a b(x_1x_2x_3x_4+ \tfrac{c}{b} x_1x_2x_5x_6+\tfrac{c}{a} x_3x_4x_5x_6).$$  This shows that $\omega^2$ is $2ab$ times a primitive class in $H^4(\zz^6)$. This is  a property that is invariant under the action of $\Aut(\zz^6)$ on $H^4 (\zz^6)$, and so $2ab$ is uniquely determined. Since $a$ is determined  by $\omega$,  and $a,b\ge 0$,  $b$ is uniquely determined.

Continuing, $\omega^3= 6 abc x_1x_2x_3x_4x_5x_6$. As before this  shows the product $6abc$ is determined  up to sign (since $c$ can be negative) and hence also $|c|$ is determined. The sign of $c$ is also uniquely determined. This is because any determinant 1 automorphism of $\zz^6$  fixes the product $6abc$ and we require that $a,b\ge 0$. \end{proof}

Thus to identify $\calg(\zz^6,\alpha)$ it suffices to consider $\alpha\in H_4(\zz^6)$ with $\alpha=\omega\cap [T]$ and $$\omega=a x_1x_2 +  bx_3x_4 +  c x_5 x_6$$ where $a|b|c$ and $a,b\ge 0$.  Here are the results.

\medskip

\begin{theorem}\label{nequals6} Given $\alpha\in H_4(\zz^6)$, there exists unique integers $a,b,c$ satisfying $a|b|c$ with $a,b\ge 0$  and $A\in \Aut(\zz^6)$ so that $A_*(\alpha)=(ax_1x_2+b x_3x_4+ c x_5x_6)\cap [T]$.  Write $q_{a,b,c}(\sigma)$ for $q_{\zz^6,   \alpha }(\sigma)$.

Then $q_{a,b,-c}(\sigma)= q_{a,b,c}(-\sigma)$ and the following equalities and estimates hold: \begin{enumerate} \item  $q_{0,0,0}(\sigma)=20 +|\sigma|$. \item $q_{1,0,0}(\sigma)=14 +|\sigma|$. \item  $q_{1,1,0}(\sigma)=10 +|\sigma|$. \item \begin{equation*}q_{1,1,1}(\sigma) = \begin{cases} 6& \text{ if } \sigma=-2 \text{ or } 0\\ 7&\text{ if }\sigma=-1 \text{ or } 1\\ 4-\sigma & \text{ if } \sigma\leq -2\\ r + \sigma& \text{ if } \sigma\ge 2, \text{ for some  }r\in \{4,6\} \end{cases} \end{equation*}

\item For general $a,b,c$, $$ \max\{s, 4+|\sigma|\}\leq q_{a,b,c}(\sigma)\leq s + |\sigma|$$ where $s=20$ if $a>1$, $s=14$ if $a=1$ and $b>1$,   and $s=10$ if $a=b=1$ and $c>1$. \end{enumerate}

\end{theorem}

\begin{proof}  The first assertion is just a restatement of Theorem \ref{6normalform}.  Replacing $(M,f)$ by $(-M,f)$ replaces $\alpha$ by $-\alpha$, and hence  $q_{a,b,c}(\sigma)=q_{-a,-b,-c}(-\sigma)$.  Since we can change the signs of any two of $a,b,c$ by a determinant 1 automorphism of $\zz^6$, we see that $q_{-a,-b,-c}(-\sigma)=q_{a,b,-c}(-\sigma)$.  Thus we assume $\alpha=\omega\cap[T]$ where $\omega=ax_1x_2+bx_3x_4+cx_5x_6$ with $a,b,c\ge 0$ and $a|b|c$.

We next establish some  lower bounds. First recall that in Section~\ref{sub6} we showed that for any $\alpha\in H_4(\zz^n)$, $q_{\zz^6,\alpha}(\sigma)\ge 2-12+16=6$.

Given $(M,f)\in \calm(\zz^6,\omega)$, the 7--dimensional subspace of $I(M,f)$, \begin{equation*}\label{veeone}V_1=\langle x_1x_2,x_1x_3,x_1x_4,x_1x_5,x_3x_5, x_3x_6, x_4x_5\rangle\end{equation*} is isotropic.   In fact, the cup product of any two elements in $V_1$ with $\omega$ vanishes. Thus  the intersection form of $M$ has a 7--dimensional isotropic subspace and hence $\beta_2(M)\ge 14 +|\sigma(M)|$. Thus \begin{equation}\label{eq6.11}q_{a,b,c}(\sigma)\ge \max\{6, 4+|\sigma|\} \end{equation} for any $a,b,c$.  (This argument offers a different perspective on the proof of Proposition~\ref{sym}, setting up needed notation.)

Now suppose that $c=0$. Given $(M,f)\in \calm(\zz^6,ax_1x_2+ bx_3x_4)$, the 10--dimensional subspace of $I(M,f)$, 
$$V_2
=\langle x_1x_2,x_1x_3,x_1x_4,x_1x_5,x_1x_6, x_2x_3,x_2x_4, x_3x_4,x_3x_5,x_3x_6\rangle, $$ 
 is isotropic.  Thus $$q_{a,b,0}(\sigma)\ge 10+|\sigma|.$$

One can get a bit more information from the space $V_2$,  Suppose that $c>1$. Then if $p$ is a prime  which divides $c$, working mod $p$ one sees that the $\zz/p$ reduction of $V_2$ is isotropic. This allows us to conclude $$q_{a,b,c}(\sigma)\ge 10 \text{ if } |c|>1.$$

Now consider the case $b,c=0$.  Given $(M,f)\in \calm(\zz^6,ax_1x_2)$, the 12-dimensional subspace of $I(M,f)$, $$V_3=\langle x_1x_2,x_1x_3,x_1x_4,x_1x_5,x_1x_6, x_2x_3,x_2x_4,x_2x_5,x_2x_6,x_3x_4,x_3x_5,x_3x_6\rangle$$ is isotropic. Hence $$q_{a,0,0}(\sigma)\ge  14+|\sigma|.$$

More generally, if $b>1$, working mod $p$ for a prime $p$ dividing $b$ one finds that $V_3$ is isotropic, so that $$q_{a,b,c}(\sigma)\ge 14 \text{ if } b>1.$$

Finally, consider the case $a=b=c=0$. For any $(M,f)\in \calm(\zz^6,0)$, the intersection form vanishes on $I(M,f)\subset H^2(M)$, which has rank $15$. This implies that $$q_{0,0,0}(\sigma)\ge |\sigma|+20.$$ As before,  (or applying Theorem \ref{multiple}) one sees $$q_{a,b,c}(\sigma)\ge 20 \text{ if } a>1.$$

\medskip

Finding good upper bounds is more involved.  We begin with a general construction which covers many cases.

Start with $M'=(F_2\times F_1)\# T^4$, where $F_g$ denotes the closed oriented surface of genus $g$. Label the generators of $\pi_1(M')$ as follows: $a_1,a_2,a_3,a_4$ generate $\pi_1(F_2)$, $b_1,b_2$ generate $\pi_1(F_1)$, and $c_1,c_2,c_3,c_4$ generate $\pi_1(T^4)$. Thus the $a_i$ commute with the $b_i$, the $b_i$ commute, the $c_i$ commute,  and $[a_1,a_2][a_3,a_4]=1$.

For $a,b,c$ as above, choose a map $f'\co M'\to T^6$ which induces the map $$a_1\mapsto \gamma_1,a_2\mapsto \gamma_2^{b},a_3\mapsto \gamma_3, a_4\mapsto \gamma_4^a,$$ $$b_1\mapsto \gamma_5,b_2\mapsto \gamma_6,$$ $$c_1\mapsto \gamma_1^{c},c_2\mapsto \gamma_2,c_3\mapsto \gamma_3, c_4\mapsto \gamma_4$$ on fundamental groups.  Then $f'_*([M'])=(a x_1x_2 + b x_3x_4 + c x_5 x_6)\cap [T^6]$. Note that $\chi(M')=-2$ and $\sigma(M')=0$.  This is true for any triple $a,b,c$ of integers; regardless of divisibility.

Next perform four surgeries on $M'$ to kill the elements $a_1^{c} c_1^{-1}$, $a_2 c_2^{-b}$, $ a_3 c_3^{-1}$, and $a_4c_4^{-a}$.  The map to $T^6$ extends over the resulting manifold which we denote by $M''$. The fundamental group of $M''$ is generated by $a_1, c_2, c_3, c_4, b_1,$ and $b_2$, which are mapped to $\gamma_1,\gamma_2,\cdots, \gamma_6$ respectively. Note that $\chi(M'')=6$ and $\sigma(M'')=0$.

For general $a,b,c$, the commutators $$\hskip-.8in z_1=[a_1,c_2] , z_2=[a_1,c_3], z_3=[a_1, c_4], z_4=[c_2,b_1],$$ $$\hskip.8in z_5=[c_2,b_2], z_6=[c_4,b_1], z_7=[c_4,b_2],$$ in $\pi_1(M'')$ need not be trivial (but $[c_3,b_1]=1=[c_3,b_2]$).  Six surgeries on $M''$ to kill  these 6 commutators yields $(M,f)\in \calm(\zz^6,\omega)$ with $\chi(M)=20$ and $\sigma(M)=0$, for any $a,b,c$. Thus we see that $q_{a,b,c}(\sigma)\leq 20 +|\sigma|$ for any $a,b,c$.  Combined with the lower bounds derived above we conclude that $$q_{0,0,0}(\sigma)=20 + |\sigma|$$ and $$q_{a,b,c}(0)=20\text{ for }a>1.$$

To treat the case $(a,b,c)=(1,b,c)$, it is convenient instead to realize $\omega=a x_1x_2 + bx_3x_4 + x_5x_6$, and then to change coordinates.  Suppose, then, that $c=1$ and $a,b$ are arbitrary. Then $a_1=c_1$ in $\pi_1(M'')$.  Thus $z_1=z_2=z_3=1$ in $\pi_1(M'')$. Hence only four surgeries are required to abelianize   $\pi_1(M'')$, yielding $(M,f)\in \calm(\zz^6,\omega)$ with $\chi(M)=14$ and $\sigma(M)=0$.  One can find a determinant $1$ automorphism that takes $\omega=a x_1x_2 + bx_3x_4 + x_5 x_6$ to $x_1x_2 + bx_3x_4 + ax_5x_6$. Thus $q_{1,b,c}(\sigma)\leq 14 + |\sigma|$.  Combined with the lower bounds derived above we conclude that $$q_{1,0,0}(\sigma)=14+|\sigma|$$ and $$14 \leq q_{1,b,c}(\sigma)\leq 14 + |\sigma|$$ when $b>1$. We will improve this below when $b=0$.

When $a=1=b$, then in $\pi_1(M'')$, $[a_3,a_4]=[c_3,c_4]=1$, and so $1=[a_1,a_2]=[a_1,c_2]=z_1$. Moreover, $z_4=[c_2,b_1]=[a_2,b_1]=1$ and similarly $z_5=1$.  Clearly $z_6=z_7=1$.  Thus to abelianize $\pi_1(M'')$ requires only two surgeries, yielding $(M,f)$ with $\chi(M)=10$ and $\sigma(M)=0$. This implies that $q_{1,1,0}(\sigma)\leq 10+|\sigma|$ and so combined with the lower bounds one concludes $$q_{1,1,0}(\sigma)= 10 + |\sigma|$$ and $$q_{1,1,c}(0)=10 \text{ if } c>1.$$

When $a=b=c=1$, then $\pi_1(M'')\cong \zz^6$, and so $q_{1,1,1}(\sigma)\leq 6+|\sigma|$ and $q_{1,1,1}(0)=6$. Moreover,  the manifold $S_6$ has $\chi(S_6)=6$ and $\sigma(S_6)=-2$. It follows from the bounds derived above that if $f\co S_6\to B_{\zz^6}$ induces an isomorphism on fundamental groups, $f_*([S_6])$ corresponds to $a=1, b=1, $ and $c=\pm1$. That  $c$ is in fact 1  follows from a symmetry argument, as follows. Notice that there is an automorphism of the genus three surface $F_3$ inducing the map on homology carrying the ordered set of generators of $H_1(F_3)$, $(x_1, x_2, \ldots, x_6)$ to $(x_3, x_4, x_5, x_6, x_1, x_2)$.  This induces an orientation preserving homeomorphism of $S_6$ inducing a similar map on $\pi_1$.  Thus, precomposing the map $S_6 \to T^6$ with this homeomorphism does not affect which class is represented in $H_4(T^6)$, but induces a map that carries $ax_1x_2 + bx_3x_4 +cx_5x_6$ to $bx_1x_2 + cx_3x_4 +ax_5x_6$.  Hence, $a = b = c$.    

We now have that  $q_{1,1,1}(\sigma)\leq 6 + |\sigma+2|$.   Thus $q_{1,1,1}(\sigma)\leq 6 + |\sigma|$ and $q_{1,1,1}(-2)=6$. These estimates are assembled in the equation \begin{equation*}q_{1,1,1}(\sigma) = \begin{cases} 6& \text{ if } \sigma=-2 \text{ or } 0\\ 7&\text{ if }\sigma=-1 \text{ or } 1\\ 4-\sigma & \text{ if } \sigma\leq -2\\ r + \sigma& \text{ if } \sigma\ge 2, \text{ for some  }r\in \{4,6\} \end{cases} \end{equation*}

\medskip

When $a=1, b=c=0$ (that is, $\omega=x_1x_2$), consider $M'=T^4\# S^1\times S^3\# S^1\times S^3$ with fundamental group generated by $a_1,a_2,a_3,a_4, b,c$.  Choose a map $f'\co M'\to T^6$ which induces the map $$a_1\mapsto \gamma_3, a_2\mapsto \gamma_4, a_3\mapsto \gamma_5, a_4\mapsto \gamma_6, b\mapsto \gamma_1,c\mapsto \gamma_6.$$ Then it is straightforward to see that $f'_*([M'])=x_1x_2\cap [T]$, that $\chi(M')=-4$, and $\sigma(M')=0$.  Nine surgeries are required to abelianize the fundamental group of $M'$; each surgery increases $\chi$ by $2$. This yields $(M,f)\in \calm(\zz^6,x_1x_2)$ with $\chi(M)=14$ and $\sigma(M)=0$.  Hence $q_{1,0,0}(\sigma)\leq |\sigma| + 14$. Combining this with the lower bound yields $$q_{1,0,0}(\sigma)= |\sigma|+14.$$

\end{proof}

One consequence of the analysis we carried out for the group $G=\zz^6$ is that the behavior of $q_G(\sigma)$ and $q_{G, \alpha}(\sigma)$ is different than when $G=\zz^n, n< 6$ or when $G$ a 2-- or 3--manifold group, since in those cases $q_G(\sigma)=k+|\sigma|$ whereas     $q_{\zz^6}$ has more than one local minimum; in fact $q_{\zz^6}(s)<q_{\zz^6}(s\pm1)$ for $s=-2,0,2$.   We will investigate this further in the next section.

\medskip
 
Calculations and estimates of $q_{\zz^n}(\sigma)$ for larger $n$ are possible, but become more difficult. In particular we have not suceeded in finding a normal form for $\omega\in H^{n-4}(\zz^n)$ for $n>6$. However, using the calculations above, the manifolds constructed in \cite{kl},  the $S_{2k}$ of Proposition \ref{sym},   Theorem \ref{exterior},  and the 7-dimensional isotropic subspace of Equation (\ref{5iso}), the following calculations are straightforward. We omit the explanations.

\begin{enumerate}
\item  $p(\zz^n)= 2,0,0,2,0,6, 4,2$ for $n=0,1,2,3,4,5,6,7$ respectively. 
\item $0\leq p(\zz^n)\leq 2-2n + C(n,2) + \epsilon_n= \frac{n^2}{2} - \frac{5n}{2} + 2 + \epsilon_n$ for all $n$, 
\item $p(\zz^{n})\leq4 - \frac{5n}{2} + C(n,2)= \frac{n^2}{2}  -\frac{7n}{2} + 4$ for $n$ even. 
\end{enumerate}

\section{Bounds on $p(\zz^n)$} \label{asym}

When we began our work in~\cite{kl} one of the first objectives was to answer a question asked by Weinberger: is $q(\zz^n)$ asymptotic to $n^2$ or $n^2 / 2$. More precisely, it followed quickly from~\cite{hw} that $$\frac{n^2 - 5n +4}{2} \le q(\zz^n) \le n^2 - 3n +2.$$  The main result of~\cite{kl} is that $q(\zz^n)$ is always within 1 of the lower bound for $n \ge 6$.  Thus, $\lim_{n\to \infty}q(\zz^n)/n^2=\tfrac{1}{2}$.

Here we would want to consider the similar question for $p(\zz^n)$. Basic estimates   show that for $n\ge 6$, $$0 \le p(\zz^n) \le  \frac{n^2 - 5n +6}{2}.$$  The lower bound comes from the fact that for a manifold $M$ with $\pi_1(M) \cong \zz^n, \sigma(M) \le \beta_2(M) - 2(n-1)$, since, by  Theorem \ref{exterior},  $H^2(M)$ contains an isotropic subspace of dimension $n-1$.  The upper bound comes about from the fact that the manifolds constructed in~\cite{kl} to realize the lower bounds of $q(\zz^n)$ have signature 0.

The use of the manifold $-S_{2k}$ introduced in Proposition \ref{sym} lets us improve the estimate for $p(\zz^n)$ in the case of $n$ even: $$ 0 \le p(\zz^n) \le \frac{n^2 - 6n + 8}{2}.$$  Even with this improvement we see that the upper bound for $p(\zz^n)$ remains asymptotic to $n^2/2$ in the sense that each inequality implies that  as $n$ goes to infinity, the quotient of the upper bound and $n^2$ is 1/2.  We cannot identify the actual behavior of $p(\zz^n)$ for large $n$, but can show that it is not asymptotic to the upper bound.

\begin{theorem}\label{13/28thm}  For an infinite set of $n$, $p(\zz^n)/n^2 \le  \frac{13}{28}$. \end{theorem}

The proof of this result is a delicate construction, efficiently building examples for large $n$  from connected sums of symmetric products of surfaces. The combinatorics is surprisingly best described in terms of projective planes over finite fields, and we begin with a detour into some of the structure of finite projective spaces.

\subsection{Points and lines in finite projective space} Let $\ff$ denote the finite field with  $p$ elements, where $p$ is a prime power.   (Our best estimate of $p(\zz^n)$ will come from setting $p=7$.)   Let $P^k$ denote $k$--dimensional projective space over $\ff$.

\begin{theorem}$\ $

\begin{enumerate} \item The number of points in $P^k$ is $n= \frac{ p^{k+1} -1}{p-1}$. \item The number of lines in $P^k$ is $L = \frac{(p^{k+1} -1)(p^{k} -1) }{ (p+1)(p-1)^2}$. \end{enumerate} \end{theorem}

\begin{proof} (1)  Recall that $P^k$ is the quotient of $\ff^{k+1} - 0$ by a free action of $\ff - 0$.  This gives the number of points in $P^k$.  Notice that if $k =1$, the projective line $P^1$ contains $p +1$ points.

(2)  Denote the number of points in $P^k$ by $n$.  Each pair of distinct points in $P^k$ determines a unique line, and thus we get $C(n,2)$ lines.  This count has repetitions, with each line being counted once for each pair of distinct points in that line, and each line has  $p+1$ points in it.  Thus, the number of lines is $C(n,2) / C(p+1,2)$. Expanding, this can be written as:

$$ \frac {  (\frac{p^{k+1} -1}{p-1} )(\frac{p^{k+1} -1}{p-1} - 1   )   } {(p+1)(p)}$$ Algebraic simplification gives the desired result.

\end{proof}

\subsection{Building examples} Fix a prime power  $p$     and integer $k$, let $n$ be the number of points in $P^k$ and let $L$ be the number of lines in $P^k$.  Let $X$ be a 4--manifold with $\pi_1(X) \cong \zz^{p+1}$, $\beta_2(X) = \beta_2$, $\chi(X) = \chi$,  and $\sigma(X) = \sigma$.

For each line $l \subset P^k$, let $X_l$ denote a copy of $X$. We will index the generators of $\pi_1(X_l)$ with the points of $l$, in arbitrary order.   Let $Y$ be the connected sum of the all the $X_l$.

Now, perform surgeries on $Y$ to identify certain pairs of generators:  if a generator of $X_{l_1}$ and a generator of $X_{l_2}$ are indexed with the same point of $P^k$, use a surgery to identify these.   Call the resulting manifold $W$.

\begin{theorem}\label{fields} $\pi_1(W) \cong \zz^n$, $ \beta_2(W) =  L\beta_2$, and $\sigma(W) = L\sigma$. \end{theorem}

\begin{proof} Clearly $\pi_1(W)$ has $n$ generators, and the first homology is $\zz^n$.  We claim that $\pi_1(W)$ is abelian.  Given any two generators, they are indexed by two points in $P^k$.  These two points determine a line $l$ in $P^k$.  Thus, the generators commute, because they both had representatives on $X_l$, which had abelian fundamental group.

Since signature adds under connected sum and surgery doesn't change the signature, we have $\sigma(W) = L\sigma(X)$, as desired.

Since $Y$ is a connected sum,  $\beta_2(Y) = L \beta_2(X)$.  Each surgery is along a curve that is of infinite order in homology, so the surgeries do not change the second Betti number:  $\beta_2(W) = L\beta_2$.

\end{proof}

Our goal is to attain asymptotics for an upper bound for $p(\zz^n)$. As our manifold $X$ we use 
$X = -S_{ p+1} $ (see Proposition \ref{sym})  so that $\pi_1(X) = \zz^{p+1}$, $\beta_2(X) =  (p^2+p+2 )/2$, $\sigma(X) = (p-1)/2$.   Noting that $p(\zz^n) \le \chi(W) - \sigma(W) \le \beta_2(W) - \sigma(W)$, we compute $$\beta_2(W) - \sigma(W)=  L(\frac{p^2 +3}{2}).$$

Dividing by $n^2$ this upper bound simplifies to be $\frac{p^2 + 3}{2(p^2 +p)}$ when terms that go to zero as $k$ increases are removed.  An elementary calculus exercise applies to show that the minimum of this function among primes occurs at either $p=5$ or $p=7$, and then a calculation shows the minimum occurs at $p=7$, where the limit is $\frac{13}{28}$.  This completes the proof of Therorem \ref{13/28thm}.\qed

\subsection{Local minima for  $q_{\zz^n}(\sigma)$}

Let 
$$C=\{ (k-\ell,k+\ell )\ | \ k,\ell\text{ non-negative integers } \}.$$
Thus $C$ consists of all integer lattice points in the $x$-$y$ plane whose coordinates have the same parity and which satisfy $ y \ge |x|$.

From Corollary \ref{basiccor}  we know that $\calg(G)$ is a union of cones of the type 
$$C_{a,b}=\{(a,b) \} + C$$ 
since if $(a,b)\in \calg(G)$, so is $(a,b)+(c,d)$ for any $(c,d)\in C$.

\begin{definition} A positive integer $a$ is called an {\em minimum point} of $q_{G} $  if  $q_{G} $ has a local minimum at $a$, i.e.~$q_G(a+1)=q_G(a-1)=q_G(a)+1$.  \end{definition}

The minimum points determine the geography completely, and correspond to irreducible manifolds when $G$ is not a free product, as the following theorem shows.

\begin{theorem}\label{except} For any $G$, $q_G $ has finitely many minimum points. Moreover, $\calg(G)$ is the union of the cones $C_{a,q_G(a)}$ as $a$ runs over the minimum points.  If a 4-manifold $M\in \calm(G)$ represents a  minimum point $a$, that is, $(\sigma(M),\chi(M))=(a, q_G(a))$, then  if $M$ is  homotopy equivalent  to a connected sum $N\# X$ for a  4-manifold $N$ and simply connected 4-manifold $X$ then $X$ is homeomorphic to the 4-sphere. In particular, if $G$ is not a free product, $M$ is irreducible as a topological 4-manifold.  
\end{theorem}
\begin{proof} The function $q_G(\sigma)-\sigma$ is integer valued  and decreasing.  As $\sigma\to\infty$, Theorem \ref{qthm}, Part 5  implies that     $q_G(\sigma)-\sigma$ is bounded below by $p(G)$. Thus there are finitely many $a\ge 0$ so that $q_G(a)-a <q_G(a-1)-(a-1)$.  Since any minimum point $a$ with $a>0$ satisfies $q_G(a-1)=q_G(a)+1$, it follows that there are finitely many minimum points $a$ with $a\ge 0$.  A similar argument using $q_G(\sigma)+\sigma$ shows that there are finitely many negative minimum points.

If $(x,q_G(x))$ is in a $C_{a,q_G(a)}$ for some minimum point $a$, then so is $(x,y)$ for any $y>q_{G}(x), y\equiv x\mod{2}$, hence to prove the second assertion it suffices to show  that for each integer $x$, $(x,q_G(x))$ lies in $C_{a,q_G(a)}$ for some minimum point $a$.   If $x$ is an minimum point the conclusion is obvious. If $q_G(x+1)<q_G(x)$,  then since $q_G$ is bounded below by $q(G)$  and since $q_G(\sigma+1)=q_G(\sigma)\pm 1$ for all $\sigma$ (Theorem \ref{qthm}, Parts 4 and 6)  there is an integer $z$ greater than $x$ so that $q_G$ has a local minimum at  $z$. If $a$ denotes the least integer greater than $x$ so that $q_G$ has a local minimum at $a$, then clearly $(x,q_G(x))\in C_{a,q_G(a)}$.  A similar argument applies if $q_G(x-1)<q_G(x)$. 

Let  $M\in \calm(G)$  satisfy $(\sigma(M),\chi(M))=(a,q_G(a))$ for some $a$.  If $M$ is homotopy equivalent to $N\# X$ where $X$ is simply connected, then $\chi(N)=\chi(M)-\beta_2(X)$ and
$\sigma(N)=\sigma(M)-\sigma(X)$. 

Let $k=\sigma(X)$.  If $k=0$, then  $\sigma(N)=\sigma(X)$, and so  $\chi(N)\ge \chi(M)$, since $\chi(M)=q_G(a)$. Thus $\beta_2(X)=0$ and so $X$ is a homotopy 4-sphere, and so by  the 4-dimensional topological Poincar\'e conjecture  (\cite{fq}), $X$ is homeomorphic to $S^4$. 

If $k>0$, then $N\#_{k-1}\cc P^2$ has signature equal to $\sigma(M)-1$ and Euler characteristic equal to $\chi(M)-\beta_2(X) + (k-1)\leq \chi(M)-1=q_G(a)-1$.  Thus $q_G(a-1)\leq q_G(a)-1$ and so $a$ is not a minimum point. A similar argument using $-\cc P^2$ shows that if $k<0$ then $a$ is not a minimum point.

\end{proof}

One can make the same definition of minimum points for $q_{G,\alpha}$ for any $\alpha\in H_4(G)$.  The assertions of Theorem  \ref{except} extend with the same proofs.  Note that for $q_G$ (but not necessarily $q_{G,\alpha}$) the set of minimum points are symmetric with respect to $0$, as one sees by reversing orientation.  Clearly the main challenge of understanding  the geography problem for a group $G$ is to identify all the minimum points (and finding the corresponding manifolds) of $q_{G}$ (and  of $q_{G,\alpha}$).

For example, for $G$ any 2-- or 3--manifold group, or  $\zz^n$ when $n<6$, $a=0$ is the only minimum point.     For $G=\zz^6$, the minimum points are exactly $-2,0, $ and $2$.  For any  $\zz^n$, $0$ is a  minimum point. 

 The construction of Theorem \ref{fields} gives groups with  at least 5 minimum points. For example, taking $G=\zz^{156}$, the example from \cite{kl} has $\sigma=0$ and $\chi=11780$. The manifold $-S_{156}$ of Proposition \ref{sym} has $\sigma=77$ and $\chi=11781$. The construction of Theorem \ref{fields} taking $p=5$ and $k=3$ yields a manifold with fundamental group $\zz^{156}$, $\sigma=1612$ and $\chi=12586$. It follows straightforwardly that $q_{\zz^{156}} $ has at least two positive minimum points ($a=77$ is a minimum point,  and there is at least one minimum point $a$ satisfying $79\leq a\leq 2408$), and from symmetry that $q_{\zz^{156}} $ has at least  five minimum points.  More generally,  it can be shown that for any $N$ there exists an $n$ so that $q_{\zz^n}$ has at least $N$ minimum points.


\section{Problems}\label{problemsec} 

\begin{enumerate}

\item  Are there groups for which the geography $\calg(G)$ depends   on the choice of category of spaces: smooth 4--manifolds, topological 4--manifolds, and 4--dimensional Poincar\'e duality spaces?  Similarly for $\calg(G,\alpha), \alpha \in H_4(G)$.

\item Develop techniques specific to the smooth category to analyze   $\calg(G)$.  Formulate the 4--dimensional topological surgery theory needed to relate the problems of computing $\calg(G)$ in the categories of topological 4--manifolds and 4-dimensional Poincar\'e complexes.

\item Are there classes of groups for which the determination of $\calg(G)$ is algorithmic?

\item  For nontrivial classes $\alpha \in H_4(G)$, find relationships between $\calg(G,\alpha)$ and $\calg(G, k\alpha), k > 1$.

\item It seems unlikely that $\calg(G,\alpha)$ is always symmetric with respect to reflection through the $\chi$ axis when $\alpha$ is non-trivial. Find a counterexample. Is $\calg(\zz^6, (x_1x_2+x_3x_4+x_5x_6)  \cap [T])$ symmetric?

\item Find relations between properties of a group $G$ and the geography of $G$.  As an example, according to~\cite{luck}, if $G$ is amenable, $q(G) \ge 0$ and $p(G)\ge 0$.

\item  How is the geography of a free product, $\calg(G_1 * G_2),$ related to 
$\calg(G_1)$ and $\calg(G_2)$?  Using connected sums one has $q(G_1 *G_2) \le q(G_1) + q(G_2) -2$.   
 For instance, we have $q(\zz^6) = 6$  but do not know whether $q(\zz^6 * \zz^6) =   9, $ or $10$.  The predicted value using additivity would be 10.

\item A class $\omega \in H^{n-4}(\zz^n)$ determines an even symmetric bilinear form 
$$\phi\co  H^2(\zz^n) \times H^2(\zz^n) \to \zz $$ by $x y   \omega = \phi(x,y) [T]$ where $[T]$ is a chosen generator of $H^n(\zz^n)$.  What unimodular forms can arise in this way?  The first unresolved case is for $n=8$, where $H^2(\zz^8) \cong \zz^{28}$ and we do not know whether the form $10H \oplus E_8$ can occur.  

\item Suppose $\omega \in H^{n-4}(\zz^n)$ determines a unimodular form $\phi$. Does there exist a smooth or topological 4-manifold with intersection form $\phi$? Does there exist a 4-dimensional Poincar\'e complex with intersection form $\phi$?

\item Determine asymptotic behavior of the geography of $\zz^n$ in terms of $n$.  For instance, the main result of~\cite{kl} is that $q(\zz^n)$ is roughly asymptotic to $ {n^2}/{2}$.  It follows that $p(\zz^n)$ is also bounded above (asymptotically) by   ${n^2}/{2}$, but according to Theorem~\ref{13/28thm} this is not the best possible: there are arbitrarily large $n$ for which $p(\zz^n) \le \frac{13}{28}n^2$.  Is it possible that $p(\zz^n)/n^2$ goes to 0 for large $n$?  

\item Determine $\calg(\zz^n, \alpha)$ for some $n$ and nonprimitive class $\alpha \in H_4(\zz^n)$.  The first case is $\calg(\zz^4, 2[T])$ where $[T]$ is a generator of $H_4(\zz^4) \cong \zz$.

\end{enumerate}


\newcommand{\etalchar}[1]{$^{#1}$} 

\end{document}